\documentclass[a4paper,fleqn]{article}
\usepackage{tikz}
\usepackage{url}
\usepackage{float}
\usepackage{placeins}
\usepackage{changepage}
\usepackage{booktabs}
\usepackage{algorithm}
\usepackage{algpseudocode}
\usepackage{hyperref}
\usetikzlibrary{arrows.meta, positioning, shapes.misc, shadows.blur}
\usepackage{graphicx}
\usepackage{caption}
\usepackage{subcaption}
\usepackage{pgfplots}
\usepackage{pgfplotstable}
\usepackage{amsmath,amssymb}
\usepackage{siunitx}
\usepackage{xcolor}
\usepackage[utf8]{inputenc}
\usepackage[T1]{fontenc}
\usepackage{lmodern}
\usepackage{geometry}
\usepackage{setspace}
\usepackage{authblk}
\usepackage[numbers,sort&compress]{natbib}
\usepackage{array}
\usepackage{longtable}
\usepackage{multirow}
\pgfplotsset{compat=1.18}
\graphicspath{{elip_paper_figures/}}

\title{\vspace{-2cm}\textbf{Approximating the Ellipse Perimeter by ``Bending the Line'': The ELIP Closed-Form Family}}
\author[1]{Salvador E. Ayala-Raggi\thanks{Corresponding author: \texttt{salvador.raggi@correo.buap.mx}}}
\author[1]{Manuel Rendón-Marín\thanks{\texttt{manuel.rendon@correo.buap.mx}}}
\affil[1]{Facultad de Ciencias de la Electrónica, Benemérita Universidad Autónoma de Puebla (BUAP), Puebla, México}
\date{October 2026}

\newcommand{\keywords}[1]{\vspace{0.5em}\noindent\textbf{Keywords: }#1\vspace{1em}}
\newcommand{\ELIP}{\mathrm{ELIP}}
\newcommand{\ppm}{\,\mathrm{ppm}}

\begin{document}
\vspace{-1cm}
\maketitle

\begin{abstract}
The perimeter of an ellipse is exactly represented by the complete elliptic integral of the second kind and therefore has no finite expression in elementary functions. This has motivated a long search for compact closed-form approximations that combine low computational cost with high accuracy. We develop a unified family called ELIP (Ellipse Linear Interpolated Perimeter), obtained by starting from the straight line that joins the two exact normalized endpoint values of the problem and then ``bending'' that line with analytically structured correction terms. Two complementary ELIP branches are studied. The first uses logarithmic--power correction kernels with optional deformation and exponential modulators. The second uses a rational residual together with one, two, or three logarithmic terms fixed by the exact elongated-ellipse asymptotics. A notation is introduced that explicitly separates Taylor-constrained coefficients from numerically free minimax parameters. This makes it possible to compare formulas by both accuracy and genuine fitting dimension.

Several compact members are identified. Rationalized low-complexity formulas reach maximum relative errors of about $6.84$, $3.57$, $0.904$, and $1.596\ppm$ with small algebraic/logarithmic structures. At the high-accuracy end, an $R_{5,2}+3\log$ ELIP formula with three fitted degrees of freedom has a rigorously certified maximum error below $0.001623\ppm$. A complementary construction fixes every coefficient of the rational core through exact Maclaurin matching and then compensates the resulting smooth beta-like residual by a very small multiplicative factor. For $R_{5,3}+3\log$, a one-parameter factor $1-Ah^{12}(1-h)^{7/4}$ yields a certified global bound $E_{\max}<0.000452\ppm$. The resulting formula is also exceptionally accurate near the circular limit because the rational-logarithmic core reproduces the exact perimeter series through $h^8$, while the external correction begins only at order $h^{12}$. The results show that explicit endpoint geometry, exact local series information, and the logarithmic structure of the degenerate limit can be combined into a systematic closed-form design methodology rather than treated as unrelated fitting heuristics.
\end{abstract}


\section{Introduction}\label{sec:introduction}

For an ellipse with semi-axes $a\ge b>0$, the perimeter is a classical example of a simple geometric quantity whose exact evaluation leads outside the elementary functions. The exact expression is
\begin{equation}
P_{\mathrm{exact}}(a,b)=4aE(e),\qquad E(e)=\int_0^{\pi/2}\sqrt{1-e^2\sin^2\theta}\,d\theta,
\qquad e=\sqrt{1-\frac{b^2}{a^2}}.
\label{eq:exact_perimeter}
\end{equation}
Equation~\eqref{eq:exact_perimeter} is mathematically exact, but evaluation of $E$ is more involved than a short formula containing only arithmetic operations, roots, logarithms, powers, or a few exponentials. The distinction matters whenever the perimeter must be evaluated repeatedly, for example inside optimization loops, embedded algorithms, geometric preprocessing, or engineering simulation codes \cite{AntennaElliptic,BrehlDowVAM,CompositeEllipse,HumanThermalEllipse}.

A particularly useful shape variable is
\begin{equation}
h=\left(\frac{a-b}{a+b}\right)^2,\qquad z=1-h=\frac{4ab}{(a+b)^2},
\label{eq:h_z_def}
\end{equation}
for which $h=0$ represents the circle and $h\rightarrow1$ the fully elongated limit. Normalizing the exact perimeter by $a+b$ gives
\begin{equation}
G_{\mathrm{exact}}(h)=\frac{P_{\mathrm{exact}}}{a+b}
=\pi\,{}_2F_1\!\left(-\frac12,-\frac12;1;h\right).
\label{eq:Gexact_hypergeom}
\end{equation}
The hypergeometric representation in Equation~\eqref{eq:Gexact_hypergeom} immediately produces the Maclaurin expansion
\begin{equation}
\begin{aligned}
G_{\mathrm{exact}}(h)=\pi\Bigg(&1+\frac{h}{4}+\frac{h^2}{64}+\frac{h^3}{256}
+\frac{25h^4}{16384}+\frac{49h^5}{65536}\\
&+\frac{441h^6}{1048576}+\frac{1089h^7}{4194304}+O(h^8)\Bigg).
\end{aligned}
\label{eq:exact_maclaurin}
\end{equation}
This series is central to the constructions developed below because it supplies exact local information at the circular endpoint.

Ramanujan proposed two celebrated approximations \cite{Ramanujan1914a,Almkvist1988}. In modern notation they are
\begin{equation}
P_{\mathrm{R1}}=\pi\left[3(a+b)-\sqrt{(3a+b)(a+3b)}\right],
\label{eq:R1}
\end{equation}
and
\begin{equation}
P_{\mathrm{R2}}=\pi(a+b)\left(1+\frac{3h}{10+\sqrt{4-3h}}\right).
\label{eq:R2}
\end{equation}
One reason for the enduring importance of Equation~\eqref{eq:R2} is that its expansion reproduces the exact normalized perimeter through the first several nontrivial powers of $h$. The agreement is sufficiently strong that many modern formulas use Ramanujan II as an already accurate base and fit only its small residual.

That strategy remains powerful, but it is not the only route. The present paper develops a complementary question: what if one begins not from Ramanujan II, but only from the two exact geometric endpoints of the normalized perimeter? At $h=0$, $G_{\mathrm{exact}}=\pi$; at $h=1$, $G_{\mathrm{exact}}=4$. The straight line joining these values is therefore
\begin{equation}
G_{\mathrm{ELIP},0}(h)=\pi+(4-\pi)h,
\label{eq:ELIP_baseline_intro}
\end{equation}
and the corresponding perimeter is $(a+b)G_{\mathrm{ELIP},0}$. Equation~\eqref{eq:ELIP_baseline_intro} is the starting point of the \emph{Ellipse Linear Interpolated Perimeter} (ELIP) family. The geometric interpretation is simple: the exact curve is obtained conceptually by bending this endpoint-exact line with a correction function.

The first ELIP branch investigated here describes that bending through logarithmic--power kernels, optionally modified by an internal deformation and by one or two exponential envelopes. This family is motivated by the fact that the exact ellipse perimeter has a logarithmic branch structure near the degenerate limit. Smooth algebraic terms alone must imitate that behavior indirectly, whereas an explicit logarithm provides the appropriate asymptotic mode at negligible structural cost.

A second ELIP branch separates the correction into two jobs. A rational function handles the smooth residual and provides efficient control near the circle, where exact Maclaurin constraints can be imposed algebraically. In parallel, one, two, or three explicit logarithmic terms reproduce the first non-analytic contributions of the elongated-limit expansion. This division of labor is important: the rational part no longer has to mimic the logarithmic singular structure, and the logarithmic part no longer has to absorb the global smooth curvature. The result is a sequence of compact rational--logarithmic formulas with a transparent relationship between exact analysis and numerical fitting.

A further feature of the ELIP construction is that the number of genuinely fitted constants can be made explicit. Some coefficients are fixed exactly by endpoint asymptotics or the Maclaurin series; others remain numerically free and are determined by minimax optimization. In the most constrained formulas, every rational coefficient is determined analytically. Their remaining global error often becomes a single smooth beta-like lobe. A small beta-shaped multiplicative factor can then redistribute this residual into a low-amplitude equioscillatory pattern using only one or two fitted parameters.

The main contribution of the paper is therefore not a single isolated approximation. It is a design framework that connects: (i) endpoint interpolation, (ii) exact Maclaurin information at the circle, (iii) exact logarithmic asymptotics at maximum elongation, (iv) low-degree rational residual models, (v) optional beta-like compensation, and (vi) rigorous interval certification for selected formulas. None of the ELIP formulas derived in this work uses either Ramanujan approximation as its mathematical base.

The remainder of the paper is organized as follows. Section~\ref{sec:history} reviews classical and modern formulas. Section~\ref{sec:methodology} presents the numerical and minimax methodology and develops the ELIP family. Section~\ref{sec:numerics} evaluates representative low-complexity and high-accuracy members, including comparisons near the circular limit. Section~\ref{sec:certification} summarizes the rigorous error certificates. Section~\ref{sec:conclusions} gives the principal conclusions.

\section{Historical Overview and Related Work}\label{sec:history}

\subsection{Classical endpoint and mean-based formulas}
The search for elementary ellipse-perimeter approximations predates modern numerical analysis. Mean-based estimates associated with Kepler's geometric intuition use either the arithmetic or geometric mean of the semi-axes \cite{Kepler1609,Sykora2005,Almkvist1988}:
\begin{equation}
P_{\mathrm{AM}}=\pi(a+b),
\label{eq:AM}
\end{equation}
\begin{equation}
P_{\mathrm{GM}}=2\pi\sqrt{ab}.
\label{eq:GM}
\end{equation}
Peano's weighted combination \cite{Kennedy1973} and Euler's quadratic-mean estimate improve these simple constructions:
\begin{equation}
P_{\mathrm{Peano}}=\pi\left(\frac{3(a+b)}{2}-\sqrt{ab}\right),
\label{eq:Peano}
\end{equation}
\begin{equation}
P_{\mathrm{Euler}}=2\pi\sqrt{\frac{a^2+b^2}{2}}.
\label{eq:EulerApprox}
\end{equation}
These formulas are useful historically because they expose the fundamental endpoint problem: a compact expression can be exact at the circle yet behave poorly when $b/a$ becomes very small.

\subsection{Euler--Ivory and the Maclaurin viewpoint}
Euler's series and Ivory's rearrangement led to the rapidly convergent expansion in the variable $h$ \cite{Euler1774,Michon2007}. In the notation of Equation~\eqref{eq:h_z_def}, the perimeter can be written as
\begin{equation}
P_{\mathrm{exact}}=\pi(a+b)\sum_{n=0}^{\infty}\binom{1/2}{n}^{\!2}h^n.
\label{eq:Ivory_series}
\end{equation}
Equation~\eqref{eq:Ivory_series} is identical to the hypergeometric representation in Equation~\eqref{eq:Gexact_hypergeom}. It gives a direct algebraic route to the exact Taylor data used later in the ELIP restrictions.

\subsection{Ramanujan and Cantrell}
Ramanujan's formulas, Equations~\eqref{eq:R1} and \eqref{eq:R2}, remain the best-known compact approximations. Ramanujan II is especially notable because its Taylor expansion agrees with the exact series through $h^4$. Cantrell \cite{Sykora2005} adds a single endpoint correction to Ramanujan II:
\begin{equation}
P_{\mathrm{C}}=\pi(a+b)\left(1+\frac{3h}{10+\sqrt{4-3h}}+\left(\frac4\pi-\frac{14}{11}\right)h^{12}\right).
\label{eq:Cantrell}
\end{equation}
The added term forces the correct fully elongated endpoint while disturbing the circle expansion only at very high order. This construction is a useful benchmark because it illustrates how an inexpensive global correction can preserve strong circular-limit behavior.

\subsection{Koshy and recent algebraic approximations}
Koshy proposed high-accuracy formulas derived from generalized means and perimeter properties \cite{Koshy2024}. One representative form writes a quarter perimeter as a generalized mean whose exponent depends on the axis ratio. Such formulas reach single-digit ppm accuracy with elementary operations, showing that purely algebraic approaches remain competitive when extreme accuracy is not required.

Moscato and Ciezak subsequently used symbolic regression together with analytical bounds and transformed variables to derive a high-accuracy approximation \cite{Moscato2024}. Because their construction is relevant to the modern closed-form literature, we reproduce its main sequence in the notation of the original source. Besides Ramanujan's second estimate $m_R=P_{\mathrm{R2}}$, they define the upper model
\begin{equation}
m_u=\pi(a+b)\left[1+\frac{(44/\pi-11)h}{10+\sqrt{4-3h}}\right],
\label{eq:Moscato_mu}
\end{equation}
and the transformed variable
\begin{equation}
h_1=\left(\frac{m_R-m_u}{m_R+m_u}\right)^2.
\label{eq:Moscato_h1}
\end{equation}
Their low-eccentricity Pad\'e component is
\begin{equation}
m_p=\pi(a+b)\frac{135168-85760h-5568h^2+3867h^3}
{135168-119552h+22208h^2-345h^3},
\label{eq:Moscato_mp}
\end{equation}
from which a second transformed variable is introduced,
\begin{equation}
h_2=\left(\frac{m_p-m_u}{m_p+m_u}\right)^2.
\label{eq:Moscato_h2}
\end{equation}
The symbolic-regression interpolation and final correction are then written in the source as
\begin{equation}
m_i=m_p\left(\frac{m_u}{m_p}\right)^{\left(h_1^2/615\right)^{(h_1/h_2)/a}},
\label{eq:Moscato_mi}
\end{equation}
\begin{equation}
m_{\mathrm{MC}}=m_i\left(1+\frac{h^{18}}{a^2-a^3}\right).
\label{eq:Moscato_mMC}
\end{equation}
Equations~\eqref{eq:Moscato_mu}--\eqref{eq:Moscato_mMC} follow the scale convention used in the published symbolic-regression construction, whose fitting and principal benchmarks set $b=1$. The method reports very small errors over the tested ranges, including strong behavior at both low and high eccentricity. It is methodologically different from ELIP because its compact structure is discovered by data-driven symbolic regression rather than built from endpoint interpolation and exact series constraints.

\subsection{Published exponential correction to Ramanujan II}
A recent published approach by Ayala-Raggi and Rendón-Marín applies a two-exponential correction to Ramanujan II \cite{Ayala2026}. A convenient representation is
\begin{equation}
P_{\mathrm{R2/2exp}}=\frac{P_{\mathrm{R2}}}{1-\sigma(h)\left(Ae^{-Bz}+Ce^{-Dz}\right)},
\label{eq:R22exp}
\end{equation}
with
\begin{equation}
\sigma(h)=\frac{1}{1+e^{-60(h-0.35)}},
\label{eq:R22exp_sigma}
\end{equation}
and a representative fitted set
\begin{equation}
A=3.37528\times10^{-4},\quad B=10.29662,\quad C=6.48093\times10^{-5},\quad D=40.89043.
\label{eq:R22exp_params}
\end{equation}
The reported full-range maximum error is about $0.57\ppm$. This formula is useful as a recent low-parameter reference, but its mathematical base remains Ramanujan II, unlike the ELIP constructions developed here.

\subsection{UNDERWOOD logarithmic--Chebyshev formulas}
Mamoun's second ellipse-perimeter paper introduces the UNDERWOOD family, whose essential idea is to represent the logarithmic branch explicitly using powers of $\xi=1-h$ multiplied by $L=\ln(1/\xi)$ and to approximate the remaining smooth component with Chebyshev polynomials \cite{Mamoun2_2026}. Define
\begin{equation}
\xi=1-h,\qquad u=2\xi-1,\qquad L=\ln\frac1\xi.
\label{eq:underwood_vars}
\end{equation}
A compact eight-constant F2 form is
\begin{equation}
\begin{aligned}
P_{\mathrm{UF2}}=\frac{a+b}{2}\Big(&7.0912285906T_0(u)-0.8666252031T_1(u)
+0.0498765412T_2(u)\\
&+0.0082178665T_3(u)+4.8744014\times10^{-4}T_4(u)\\
&+0.2463044164\xi^2L+0.1343197423\xi^3L+0.0194597499\xi^4L\Big).
\end{aligned}
\label{eq:underwood_F2}
\end{equation}
The paper reports a practical error scale near $0.015\ppm$ for this compact member, with tighter Remez-refitted coefficient sets also discussed.

For comparison with our high-accuracy plots we also use the ten-constant UNDERWOOD set
\begin{equation}
P_{\mathrm{U10}}=\frac{a+b}{2}\left[\sum_{k=0}^{4}c_kT_k(u)+\sum_{j=2}^{6}d_j\xi^jL\right],
\label{eq:underwood_U10}
\end{equation}
where
\begin{equation}
\begin{aligned}
(c_0,c_1,c_2,c_3,c_4)=&(7.088680365094343,-0.868451433296990,0.0520545244008907,\\
&0.0100440873982781,0.000857760231117241),
\end{aligned}
\label{eq:underwood_U10_c}
\end{equation}
\begin{equation}
\begin{aligned}
(d_2,d_3,d_4,d_5,d_6)=&(0.249256072270624,0.165858509549193,0.047113266952720,\\
&-0.00430961785445353,0.000653807402098442).
\end{aligned}
\label{eq:underwood_U10_d}
\end{equation}
On our independent dense grid, the coefficient set in Equations~\eqref{eq:underwood_U10}--\eqref{eq:underwood_U10_d} gives a maximum error of approximately $0.0006512\ppm$. This numerical comparison is distinct from the more general theoretical claims made for other UNDERWOOD bases and larger parameter budgets.

\subsection{Popenkov's rational--logarithmic correction}

Popenkov combines Ramanujan II with an explicit logarithmic residual
and a low-degree rational denominator~\cite{Popenkov2026}. With
$s=a+b$, the formula is

\begin{equation}
P_{\mathrm{Pop}}
=
\pi s
\left(
1+
\frac{3h}{10+\sqrt{1+3z}}
+
h^5\frac{N(z)}{D(z)}
\right),
\label{eq:Popenkov}
\end{equation}

where

\begin{equation}
N(z)
=
p_0+p_1z+p_2z^2+p_3z^3
-
\frac{z^2\ln z}{8\pi},
\label{eq:Popenkov_N}
\end{equation}

and

\begin{equation}
D(z)
=
1+\frac{8}{15}z+\lambda z^2
+\frac{427}{100}\lambda z^3.
\label{eq:Popenkov_D}
\end{equation}

The coefficients reported by Popenkov are

\begin{equation}
\begin{aligned}
p_0
&=
\frac{4}{\pi}-\frac{14}{11}
\\
&\approx
0.0005122720078899588783,
\\[4pt]
p_1
&=
\frac{19}{\pi}
-\frac{1465}{242}
+\frac{8}{15}p_0
\\
&\approx
-0.005557959034898739298,
\\[4pt]
p_2
&=
0.019232178410918235337,
\\[4pt]
p_3
&=
0.012701584343839720567.
\end{aligned}
\label{eq:Popenkov_parameters}
\end{equation}

Thus, $p_0$ and $p_1$ incorporate exact analytical information,
whereas $p_2$ and $p_3$ are given numerically in the reported
rational--logarithmic approximation. The final minimax parameter is

\begin{equation}
\lambda
=
222.6232676652460423250403081343657.
\label{eq:Popenkov_lambda}
\end{equation}

Popenkov gives a rigorous interval-arithmetic certificate

\begin{equation}
1.5851297637\times10^{-8}
\leq
\max |\varepsilon|
<
1.5852\times10^{-8},
\label{eq:Popenkov_certificate}
\end{equation}

corresponding to a maximum absolute relative error below
$0.015852$ ppm.


\subsection{Error Analysis for extreme ratios \texorpdfstring{$a/b$}{a/b}}\label{sec:extreme}
Two limiting cases provide a useful structural test. In the circular limit,
\begin{equation}
\lim_{a/b\to1}\frac{P_{\mathrm{exact}}}{a}=2\pi,
\label{eq:circle_limit}
\end{equation}
while for a fully flattened ellipse,
\begin{equation}
\lim_{a/b\to\infty}\frac{P_{\mathrm{exact}}}{a}=4.
\label{eq:flat_limit}
\end{equation}
For any approximation $P^*$ we use the signed relative error
\begin{equation}
\varepsilon(h)=\frac{P^*(h)-P_{\mathrm{exact}}(h)}{P_{\mathrm{exact}}(h)},\qquad E_{\max}=10^6\max_{0\le h\le1}|\varepsilon(h)|\ \mathrm{ppm}.
\label{eq:error_definition}
\end{equation}
Table~\ref{tab:limits} summarizes the endpoint behavior. The key structural advantage of ELIP is immediate: the baseline is exact at both endpoints, and the rational--logarithmic corrections are built so that their additional terms vanish appropriately at $h=0$ and $h=1$.

\begin{table}[H]
\centering
\footnotesize
\renewcommand{\arraystretch}{1.18}
\caption{Limiting behavior at the circle and at complete elongation. Values for UNDERWOOD refer to the displayed decimal coefficient sets; ELIP endpoint equalities are structural.}
\label{tab:limits}
\begin{tabular}{p{4.2cm}ccc}
\toprule
Method & $\lim_{a/b\to1}P^*/a$ & $\lim_{a/b\to\infty}P^*/a$ & $|\varepsilon_\infty|$ (ppm)\\
\midrule
Geometric mean & $2\pi$ & $0$ & $1{,}000{,}000$\\
Arithmetic mean & $2\pi$ & $\pi$ & $214{,}602$\\
Euler & $2\pi$ & $\sqrt2\,\pi$ & $110{,}721$\\
Peano & $2\pi$ & $3\pi/2$ & $178{,}097$\\
Ramanujan I & $2\pi$ & $(3-\sqrt3)\pi$ & $4155.03$\\
Ramanujan II & $2\pi$ & $14\pi/11$ & $402.337$\\
Cantrell & $2\pi$ & $4$ & $0$\\
R2/2exp & quasi-exact & quasi-exact & $\sim1.9\times10^{-4}$\\
UNDERWOOD F2 (displayed set) & approximate & approximate & $\sim 0.0114$\\
UNDERWOOD U10 (displayed set) & approximate & approximate & $\sim5.5\times10^{-4}$\\
Popenkov & $2\pi$ & $4$ & $0$\\
ELIP families in this work & $2\pi$ & $4$ & $0$\\
\bottomrule
\end{tabular}
\end{table}

\section{Methodology}\label{sec:methodology}

\subsection{Numerical computation of the exact elliptic integral}
The numerical reference values were evaluated with arbitrary-precision arithmetic using \texttt{mpmath}. For a chosen pair $(a,b)$, Equation~\eqref{eq:exact_perimeter} was evaluated through the complete elliptic integral of the second kind. The shape-only representation in Equation~\eqref{eq:Gexact_hypergeom} was also used independently when convenient. Working precision was increased substantially beyond the target errors for refined minimax calculations and for locating extrema. Dense-grid scans were used only for exploration and plotting; the certified results in Section~\ref{sec:certification} use interval arithmetic rather than a finite sampling argument.

The homogeneity of the perimeter reduces the numerical study to a single shape parameter. Since both $P_{\mathrm{exact}}$ and every formula considered here are homogeneous of degree one in $(a,b)$, the relative error depends only on $h$ (or equivalently $z$). This avoids artificial dependence on the scale of the ellipse.

\subsection{Minimax fitting}\label{sec:minimax}
When a formula contains genuinely free coefficients, they are selected by minimizing the maximum relative error over the physical interval. If $\boldsymbol\theta$ denotes the free parameter vector, the numerical objective is
\begin{equation}
\Phi(\boldsymbol\theta)=\max_{0\le h\le1}\left|\frac{P_{\mathrm{aprox}}(h;\boldsymbol\theta)}{P_{\mathrm{exact}}(h)}-1\right|.
\label{eq:minimax_objective}
\end{equation}
Analytical restrictions are applied before the numerical search whenever possible. Thus, if $k$ Taylor conditions determine $k$ coefficients algebraically, these coefficients are eliminated from $\boldsymbol\theta$ and recomputed exactly for every trial value of the remaining variables. This prevents the optimizer from wasting degrees of freedom on conditions that are known in closed form.

Algorithm~\ref{alg:minimax} summarizes the procedure. Differential evolution or another global search is used to explore the reduced parameter space, followed by derivative-free local refinement, for example using Powell's method \cite{Powell1981}, and high-precision extremum localization. Equioscillation of the signed error is used as a diagnostic of a well-resolved minimax fit, but a dense numerical minimax result is not called a rigorous certificate unless the independent interval procedure of Section~\ref{sec:certification} succeeds.

\begin{algorithm}[H]
\caption{Constrained minimax fitting used for ELIP formulas}\label{alg:minimax}
\begin{algorithmic}[1]
\State Choose an ELIP architecture and identify all nominal coefficients.
\State Derive endpoint and Taylor constraints symbolically.
\State Eliminate the analytically determined coefficients.
\State Construct a dense, endpoint-enriched grid in $h\in[0,1]$.
\State Evaluate the high-precision reference perimeter on the grid.
\State Define the objective $\Phi$ in Equation~\eqref{eq:minimax_objective} over only the free parameters.
\State Perform a bounded global search for a low value of $\Phi$.
\State Refine the best candidate with a local optimizer or one-dimensional Remez-like balancing when applicable.
\State Locate the dominant extrema of the signed error at high precision.
\State Check denominator positivity and endpoint limits.
\State Report the numerical maximum; if a theorem-level bound is required, run the independent interval certifier.
\end{algorithmic}
\end{algorithm}

Minor differences in the fitted parameter values may arise when the minimax optimization is rerun, due to grid resolution, stopping tolerances, initialization, and floating-point effects. Such variations generally lead to essentially the same error curve and maximum relative error.

\subsection{Constructing accurate ellipse perimeter approximations from a linear interpolation base: the ELIP family}\label{sec:ELIPbase}
\begin{figure}[H]
\centering
\includegraphics[width=0.78\linewidth]{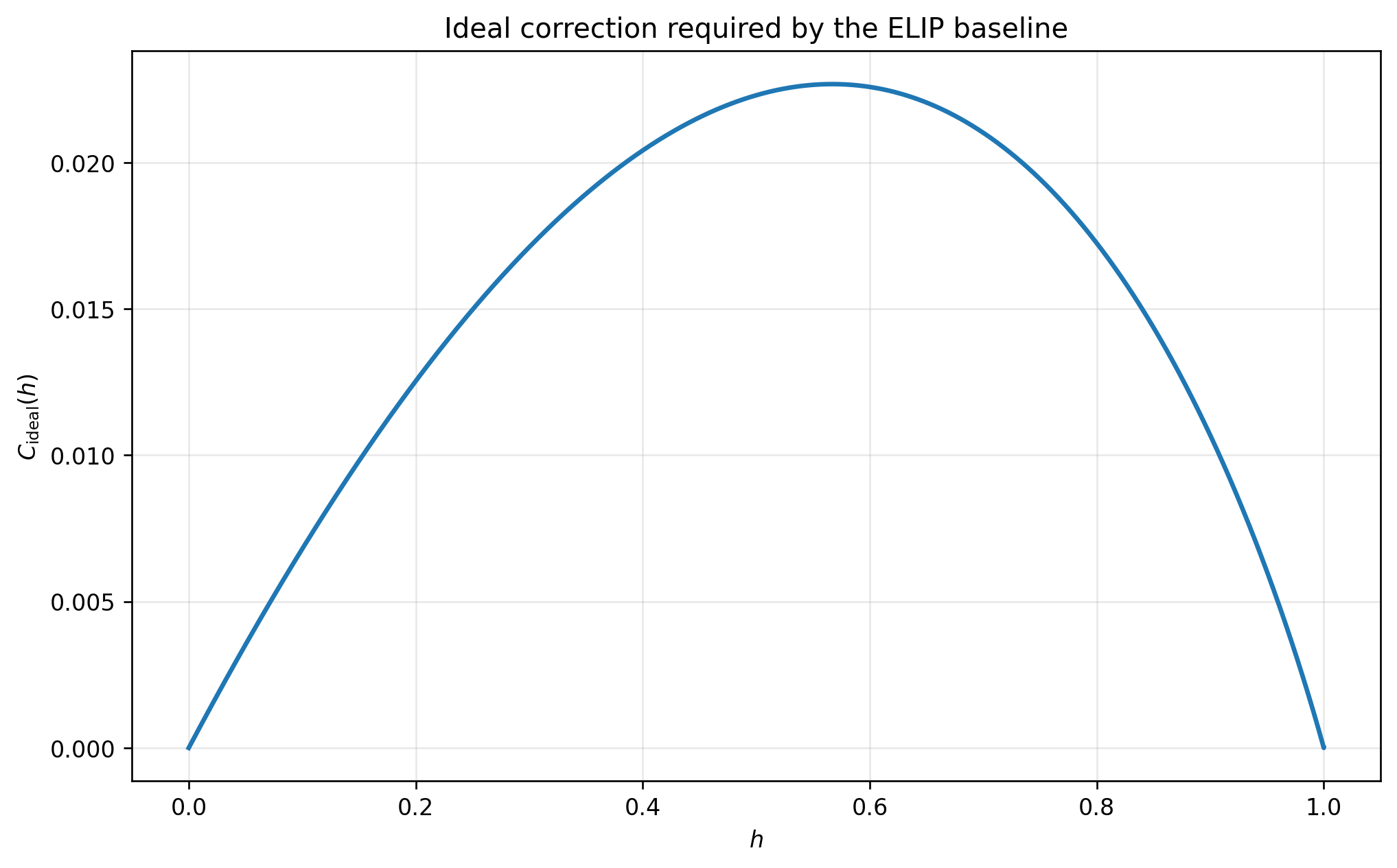}
\caption{Ideal correction required to transform the endpoint-exact ELIP line into the exact normalized perimeter.}
\label{fig:idealC}
\end{figure}

Equation~\eqref{eq:ELIP_baseline_intro} is exact at both geometric endpoints, so the exact normalized perimeter can be written as the line minus an ideal correction:
\begin{equation}
G_{\mathrm{exact}}(h)=\pi+(4-\pi)h-C_{\mathrm{ideal}}(h).
\label{eq:Cideal_definition}
\end{equation}
Therefore
\begin{equation}
C_{\mathrm{ideal}}(h)=\pi+(4-\pi)h-G_{\mathrm{exact}}(h).
\label{eq:Cideal_formula}
\end{equation}
Figure~\ref{fig:idealC} shows this correction. The phrase ``bending the line'' refers literally to replacing the straight interpolation in Equation~\eqref{eq:ELIP_baseline_intro} by Equation~\eqref{eq:Cideal_definition}, with a compact approximation to $C_{\mathrm{ideal}}$.

The generic corrected ELIP form is
\begin{equation}
P_{\mathrm{ELIP}}=(a+b)\left[\pi+(4-\pi)h-C(h)\right].
\label{eq:ELIP_generic}
\end{equation}
The challenge is therefore transferred from approximating $P$ directly to finding a compact $C(h)$ that respects both endpoints and captures the curvature of Equation~\eqref{eq:Cideal_formula}.

\subsection{Notation for ELIP architectures}\label{sec:notation}
To make structural complexity visible in the name of the formulas described in the subsequent sections, the logarithmic--power--exponential branch is denoted as
\begin{equation}
\boxed{\ELIP\;LPD E_k[T_r/F_f]},
\label{eq:notation_LPDE}
\end{equation}
where $L$ denotes the logarithm, $P$ the power envelope, $D$ an optional deformation $h\mapsto h^\rho$, and $E_k$ the number of exponential modulators. When no deformation is present the $D$ is omitted. For the special restricted six-parameter construction below, an inexpensive algebraic modulator $(1+Bh)$ is indicated by $M_1$.

On the other hand, the rational--logarithmic branch is denoted by
\begin{equation}
\boxed{\ELIP\;R_{n,m}+q\log[T_r/F_f;\beta_t]}.
\label{eq:notation_rational}
\end{equation}
Here $R_{n,m}$ is the rational degree, $q$ is the number of exact elongated-limit logarithmic terms, $T_r$ counts coefficients fixed by Taylor/Maclaurin conditions, $F_f$ counts free coefficients remaining in the rational core, and $\beta_t$ denotes an optional beta-type compensator with $t$ fitted parameters. The beta suffix is omitted when no such factor is present.

A rationalized implementation that replaces fitted decimals by small fractions is marked with the suffix ``rat''. Rationalization may slightly relax exact Taylor equalities; this distinction is stated explicitly whenever it occurs.

\subsection{Logarithmic--power and exponential ELIP corrections}\label{sec:LPDE}
The simplest useful logarithmic correction begins with
\begin{equation}
C(h)=A(1-h)^p\ln\frac1{1-h}.
\label{eq:ELIP_LP}
\end{equation}
A deformation of the internal variable gives the three-parameter architecture
\begin{equation}
C(h)=A(1-h^\rho)^p\ln\frac1{1-h},
\label{eq:ELIP_LPD}
\end{equation}
which we denote $\ELIP\;LPD[T_0/F_3]$. A minimax fit yields approximately
\begin{equation}
A=0.07280973292742082,\qquad \rho=1.023172598918835,\qquad p=1.20632845179987,
\label{eq:ELIP_LPD_params}
\end{equation}
with maximum error about $6.50\ppm$. A compact rationalized form is
\begin{equation}
P_{\mathrm{LPD,rat}}=(a+b)\left[\pi+(4-\pi)h-\frac4{55}(1-h^{41/40})^{35/29}\ln\frac1{1-h}\right],
\label{eq:ELIP_LPD_rat}
\end{equation}
whose dense-grid maximum is approximately $6.84\ppm$.

One exponential modulator without the deformation gives $\ELIP\;LPE_1[T_0/F_5]$:
\begin{equation}
C(h)=A(1-h)^p\ln\frac1{1-h}\left[1+K e^{-D(1-h)^s}\right],
\label{eq:ELIP_LPE1}
\end{equation}
with representative parameters
\begin{equation}
A=0.08035732,\quad p=1.09779322,\quad K=-0.55781159,\quad D=1.80916026,\quad s=0.41266493,
\label{eq:ELIP_LPE1_params}
\end{equation}
reaching about $0.087\ppm$. Allowing a deformation $\rho$ in the first factor produces $\ELIP\;LPDE_1[T_0/F_6]$:
\begin{equation}
C(h)=A(1-h^\rho)^p\ln\frac1{1-h}\left[1+K e^{-D(1-h)^s}\right],
\label{eq:ELIP_LPDE1}
\end{equation}
with
\begin{equation}
\begin{aligned}
A&=0.0814628277594858, & \rho&=1.000677047628661, & p&=1.09242881461921,\\
K&=-0.5860457228132984, & D&=1.730433583994363, & s&=0.4040670464083422,
\end{aligned}
\label{eq:ELIP_LPDE1_params}
\end{equation}
and a numerical maximum near $0.0292\ppm$.
Another way to improve Eq. \ref{eq:ELIP_LPE1} is by adding a simple linear factor. The correction is written as
\begin{equation}
C(h)
=
A(1-h)^p
\ln\!\left(\frac{1}{1-h}\right)
\left[
1+K e^{-D(1-h)^s}
\right]
(1+Bh),
\label{eq:restricted_six}
\end{equation}
and denoted as $\ELIP\;LPE_1M_1[T_0/F_6]$. A representative numerical minimax fit is
\begin{equation}
\begin{aligned}
A&=0.0826528127793086, &
p&=1.0885906179096696, &
K&=-0.6083034155074792,\\
D&=1.6511653658733876, &
s&=0.3967305433972328, &
B&=0.0028494542344265.
\end{aligned}
\label{eq:LPE1M1_free_params}
\end{equation}
For this unrestricted fit, the numerically observed maximum relative error is approximately
\begin{equation}
E_{\max}\simeq 0.01608833\ \mathrm{ppm}.
\label{eq:LPE1M1_free_error}
\end{equation}

Finally, a two-exponential form is $\ELIP\;LPE_2[T_0/F_8]$:
\begin{equation}
C(h)=A(1-h)^p\ln\frac1{1-h}\prod_{i=1}^{2}\left[1+K_i e^{-D_i(1-h)^{s_i}}\right],
\label{eq:ELIP_LPE2}
\end{equation}
with one fitted set
\begin{equation}
\begin{aligned}
A&=0.08212693, & p&=1.08665524,\\
K_1&=-0.60187957, & D_1&=1.72841795, & s_1&=0.39240755,\\
K_2&=-0.02710856, & D_2&=1.76765048, & s_2&=0.56737466,
\end{aligned}
\label{eq:ELIP_LPE2_params}
\end{equation}
for which the observed maximum is approximately $0.01195\ppm$. These formulas show that the logarithmic core is effective even before introducing the rational residual developed later.

\subsection{Taylor-restricted logarithmic--exponential ELIP}\label{sec:restricted_LPE}

A useful intermediate experiment is to compare two different ways of using the same
six-parameter logarithmic--exponential architecture of $\ELIP\;LPE_1M_1[T_0/F_6]$ in Eq. \ref{eq:restricted_six}. For the unrestricted fit, the numerically observed maximum relative error is approximately
\begin{equation}
E_{\max}\simeq 0.01608833\ \mathrm{ppm}.
\label{eq:LPE1M1_free_error}
\end{equation}
Thus, when the only objective was to minimize the largest error over the complete interval
$0<h<1$, the six available degrees of freedom can be used very efficiently.

However, this fit does not force the approximation to reproduce the exact Maclaurin
coefficients at the circular limit.  Although the formula is exact in value at $h=0$,
its first nonzero local error is not forced to begin at a high power of $h$.
This motivates a second construction in which part of the available numerical freedom is
exchanged for exact local information.

\paragraph{Exact Taylor conditions.}
The ideal correction is defined by
\begin{equation}
C_{\mathrm{ideal}}(h)
=
\pi+(4-\pi)h-G_{\mathrm{exact}}(h),
\label{eq:Cideal_restricted_definition}
\end{equation}
where
\begin{equation}
G_{\mathrm{exact}}(h)
=
\frac{P_{\mathrm{exact}}}{a+b}.
\label{eq:Cideal_restricted_G}
\end{equation}
Using the exact Maclaurin expansion of the ellipse perimeter gives
\begin{equation}
C_{\mathrm{ideal}}(h)
=
\left(4-\frac{5\pi}{4}\right)h
-\frac{\pi}{64}h^2
-\frac{\pi}{256}h^3
-\frac{25\pi}{16384}h^4
+O(h^5).
\label{eq:Cideal_taylor}
\end{equation}
The Taylor-restricted construction requires the correction in
Equation~\eqref{eq:restricted_six} to have exactly the same first four coefficients.

The coefficient matching can be organized without expanding the complete expression at once.
Let
\begin{equation}
z=1-h,
\qquad
L(h)=\ln\!\left(\frac{1}{1-h}\right).
\label{eq:restricted_zL}
\end{equation}
We first write
\begin{equation}
z^p(1+Bh)
=
1+w_1h+w_2h^2+w_3h^3+O(h^4),
\label{eq:restricted_wseries}
\end{equation}
where
\begin{equation}
w_1=B-p,
\label{eq:restricted_w1}
\end{equation}
\begin{equation}
w_2=\frac{p(p-1)}{2}-pB,
\label{eq:restricted_w2}
\end{equation}
\begin{equation}
w_3=
-\frac{p(p-1)(p-2)}{6}
+\frac{Bp(p-1)}{2}.
\label{eq:restricted_w3}
\end{equation}

Next, define
\begin{equation}
q=K e^{-D}.
\label{eq:restricted_q}
\end{equation}
The exponential term can be expanded as
\begin{equation}
e^{-D(1-h)^s}
=
e^{-D}
\left[
1+e_1h+e_2h^2+e_3h^3+O(h^4)
\right],
\label{eq:restricted_expansion}
\end{equation}
with
\begin{equation}
e_1=Ds,
\label{eq:restricted_e1}
\end{equation}
\begin{equation}
e_2=
\frac{D^2s^2-Ds(s-1)}{2},
\label{eq:restricted_e2}
\end{equation}
and
\begin{equation}
e_3=
\frac{Ds(s-1)(s-2)}{6}
-\frac{D^2s^2(s-1)}{2}
+\frac{D^3s^3}{6}.
\label{eq:restricted_e3}
\end{equation}
Therefore,
\begin{equation}
1+K e^{-D(1-h)^s}
=
v_0+v_1h+v_2h^2+v_3h^3+O(h^4),
\label{eq:restricted_vseries}
\end{equation}
where
\begin{equation}
v_0=1+q,
\qquad
v_1=qe_1,
\qquad
v_2=qe_2,
\qquad
v_3=qe_3.
\label{eq:restricted_vcoeff}
\end{equation}

Multiplying the two non-logarithmic factors gives
\begin{equation}
z^p(1+Bh)
\left[
1+K e^{-Dz^s}
\right]
=
H_0+H_1h+H_2h^2+H_3h^3+O(h^4),
\label{eq:restricted_Hseries}
\end{equation}
with
\begin{equation}
H_0=v_0,
\label{eq:restricted_H0}
\end{equation}
\begin{equation}
H_1=w_1v_0+v_1,
\label{eq:restricted_H1}
\end{equation}
\begin{equation}
H_2=w_2v_0+w_1v_1+v_2,
\label{eq:restricted_H2}
\end{equation}
\begin{equation}
H_3=w_3v_0+w_2v_1+w_1v_2+v_3.
\label{eq:restricted_H3}
\end{equation}

Finally,
\begin{equation}
L(h)
=
h+\frac{h^2}{2}+\frac{h^3}{3}+\frac{h^4}{4}+O(h^5).
\label{eq:restricted_logseries}
\end{equation}
Substitution into Equation~\eqref{eq:restricted_six} therefore gives
\begin{equation}
C(h)
=
c_1h+c_2h^2+c_3h^3+c_4h^4+O(h^5),
\label{eq:restricted_Cseries}
\end{equation}
where
\begin{equation}
c_1=AH_0,
\label{eq:restricted_c1}
\end{equation}
\begin{equation}
c_2=A\left(H_1+\frac{H_0}{2}\right),
\label{eq:restricted_c2}
\end{equation}
\begin{equation}
c_3=A\left(H_2+\frac{H_1}{2}+\frac{H_0}{3}\right),
\label{eq:restricted_c3}
\end{equation}
\begin{equation}
c_4=A\left(H_3+\frac{H_2}{2}+\frac{H_1}{3}+\frac{H_0}{4}\right).
\label{eq:restricted_c4}
\end{equation}

The four Taylor restrictions are now simply
\begin{equation}
c_1=4-\frac{5\pi}{4},
\label{eq:restricted_match1}
\end{equation}
\begin{equation}
c_2=-\frac{\pi}{64},
\label{eq:restricted_match2}
\end{equation}
\begin{equation}
c_3=-\frac{\pi}{256},
\label{eq:restricted_match3}
\end{equation}
\begin{equation}
c_4=-\frac{25\pi}{16384}.
\label{eq:restricted_match4}
\end{equation}
These four equations remove four effective degrees of freedom from the six-parameter
architecture.  Two degrees of freedom remain, which explains the notation
\begin{equation}
\ELIP\;LPE_1M_1[T_4/F_2].
\label{eq:LPE1M1_restricted_name}
\end{equation}

One convenient computational organization is to use $D$ and $s$ as the two outer free
variables.  For every trial pair $(D,s)$, Equations~\eqref{eq:restricted_match1}--
\eqref{eq:restricted_match4} are solved for $A$, $p$, $K$, and $B$.  The resulting
two-dimensional family already satisfies the four exact Taylor conditions by construction.
A minimax search is then performed only over $(D,s)$.  In this way, numerical optimization
does not determine all six coefficients independently; it only selects the two remaining
degrees of freedom after the exact local information has been imposed.

A representative restricted solution is

\begin{equation}
\boxed{
\begin{aligned}
P_{\mathrm{aprox}}
&=(a+b)\Bigg[
\pi+(4-\pi)h
-A(1-h)^p
\ln\left(\frac{1}{1-h}\right)
\left(1+K\exp\left[-D(1-h)^s\right]\right)
(1+Bh)
\Bigg],
\end{aligned}
}
\label{eq:ELIP_LPE1M1_T4F2}
\end{equation}

\begin{equation}
\begin{aligned}
A&=0.0835719073318, &
p&=1.083530519481, &
K&=-0.630322440201,\\
D&=1.606852416470, &
s&=0.392617460220, &
B&=0.002458526279.
\end{aligned}
\label{eq:restricted_six_params}
\end{equation}
Its numerically observed full-range maximum is approximately
\begin{equation}
E_{\max}\simeq0.13384\ \mathrm{ppm},
\label{eq:restricted_six_global_error}
\end{equation}
with the largest error near
\begin{equation}
h\simeq0.97675.
\label{eq:restricted_six_hmax}
\end{equation}

The price of the Taylor restrictions is therefore a loss of global minimax accuracy:
the unrestricted six-parameter fit in Equation~\eqref{eq:LPE1M1_free_params} has a
maximum error of about $0.01609\ \mathrm{ppm}$, whereas the restricted version has a
maximum of about $0.13384\ \mathrm{ppm}$.  The benefit appears near the circular limit.
For the unrestricted fit, the maximum error on $0<h<0.3$ is approximately
$0.01609\ \mathrm{ppm}$, while for the restricted formula it is only
\begin{equation}
\max_{0<h<0.3}|\varepsilon(h)|
\simeq0.00123\ \mathrm{ppm}.
\label{eq:restricted_near_circle}
\end{equation}
Thus, over this near-circle interval, imposing the four Taylor conditions improves the
worst-case error by approximately a factor of thirteen relative to the unrestricted fit.
For comparison, Ramanujan II reaches approximately $0.08970\ \mathrm{ppm}$ on the same
interval, so the restricted ELIP formula is about $73$ times more accurate there.

More importantly, the improvement is asymptotic rather than merely numerical.
Because Equations~\eqref{eq:restricted_match1}--\eqref{eq:restricted_match4} reproduce
the exact correction through $h^4$, the first possible discrepancy occurs at order $h^5$:
\begin{equation}
\frac{P_{\mathrm{aprox}}-P_{\mathrm{exact}}}
{P_{\mathrm{exact}}}
=
O(h^5),
\qquad h\rightarrow0.
\label{eq:restricted_order_h5}
\end{equation}
For the representative solution above, the magnitude of the residual $h^5$ coefficient
is approximately $1.9\%$ of the corresponding Ramanujan-II coefficient.  Hence, in the
asymptotic circular limit, the restricted formula is approximately $53$ times more accurate
than Ramanujan II at leading order.

\begin{figure}[H]
\centering
\includegraphics[width=0.94\linewidth]{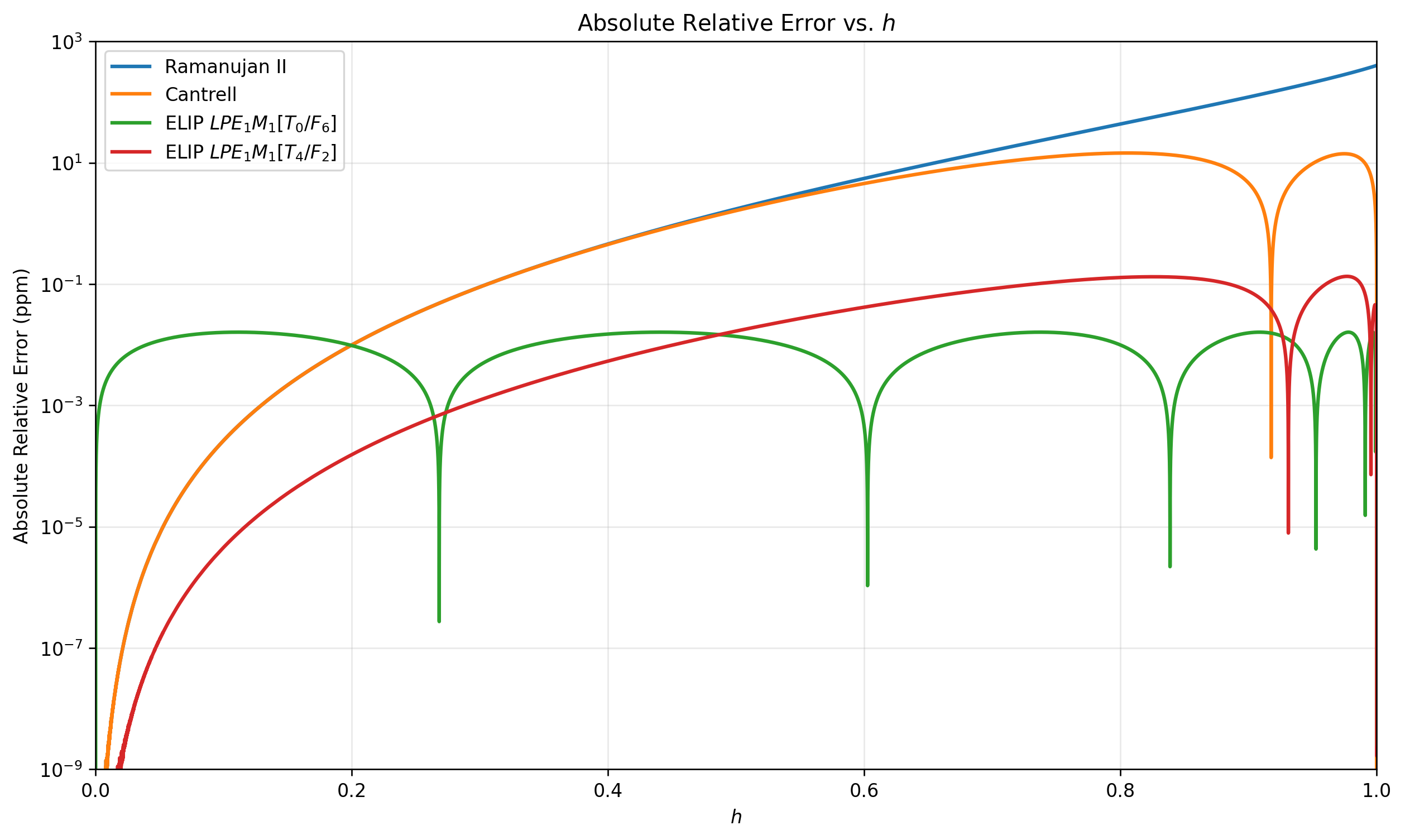}
\caption{Absolute relative error of Ramanujan II, $\ELIP\;LPE_1M_1[T_0/F_6]$, $\ELIP\;LPE_1M_1[T_4/F_2]$}
\label{fig:error_LPE1M1_T0F6_T4F2}
\end{figure}

This example illustrates the central purpose of Taylor restriction in the ELIP family.
A fully free minimax fit uses all available parameters to reduce the largest error over the
whole interval.  A Taylor-restricted fit deliberately sacrifices part of that global freedom
in exchange for exact local structure.  The result can have a larger full-range maximum while
being substantially, and asymptotically, more accurate for nearly circular ellipses. Figure \ref{fig:error_LPE1M1_T0F6_T4F2} shows the difference.

\subsection{Rational--Logarithmic ELIP Approximations for the Ellipse Perimeter}\label{sec:rationallog}

The second major ELIP branch uses
\begin{equation}
P_{\mathrm{aprox}}=(a+b)\left[\pi+(4-\pi)h-hzR_{n,m}(z)+L_q(z)\right],
\label{eq:RL_general}
\end{equation}
where
\begin{equation}
h=\left(\frac{a-b}{a+b}\right)^2,
\qquad
z=1-h,
\label{eq:RL_hz}
\end{equation}
and
\begin{equation}
R_{n,m}(z)=
\frac{\pi-3+a_1z+\cdots+a_nz^n}
{1+b_1z+\cdots+b_mz^m}.
\label{eq:Rnm_general}
\end{equation}
The constant $\pi-3$ in Equation~\eqref{eq:Rnm_general} is fixed by the elongated endpoint.
The first three exact logarithmic blocks are
\begin{equation}
L_1(z)=\frac18z^2\ln\frac1z,
\label{eq:L1}
\end{equation}
\begin{equation}
L_2(z)=
\left(\frac18z^2+\frac3{32}z^3\right)\ln\frac1z,
\label{eq:L2}
\end{equation}
\begin{equation}
L_3(z)=
\left(\frac18z^2+\frac3{32}z^3+\frac{75}{1024}z^4\right)\ln\frac1z.
\label{eq:L3}
\end{equation}
These terms come directly from the exact elongated-limit expansion of the elliptic integral.
Their analytical role is to remove, explicitly, the leading non-analytic logarithmic modes near
$z=0$, leaving a substantially smoother remainder to be represented by the rational function.

An important point in the notation used below is the meaning of quantities such as $R(1)$.
Since $z=1-h$, the circular limit $h=0$ corresponds to $z=1$. Therefore,
\begin{equation}
R(1)
=
R_{n,m}(z)\big|_{z=1}
=
\frac{\pi-3+a_1+\cdots+a_n}
{1+b_1+\cdots+b_m},
\label{eq:R_at_1_general}
\end{equation}
is simply the value of the rational correction at the circular endpoint. Conditions imposed on
$R(1)$, $R'(1)$, and higher derivatives at $z=1$ are thus equivalent to requiring that the ELIP
approximation reproduce successive coefficients of the exact Maclaurin expansion about $h=0$.

The exact normalized perimeter has the expansion
\begin{equation}
\frac{P_{\mathrm{exact}}}{a+b}
=
\pi+\frac{\pi}{4}h+\frac{\pi}{64}h^2
+\frac{\pi}{256}h^3
+\frac{25\pi}{16384}h^4+\cdots.
\label{eq:exact_circle_series_RL}
\end{equation}

For the first rational-only member, $\ELIP\;R_{2,2}[T_1/F_3]$, we write
\begin{equation}
R_{2,2}(z)
=
\frac{\pi-3+a_1z+a_2z^2}
{1+b_1z+b_2z^2}.
\label{eq:R22_structure_revised}
\end{equation}
Near the circle, $z=1-h$ and
\begin{equation}
R_{2,2}(1-h)=R_{2,2}(1)+O(h).
\label{eq:R22_local_R}
\end{equation}
Consequently,
\begin{equation}
hzR_{2,2}(z)
=
h(1-h)R_{2,2}(1-h)
=
R_{2,2}(1)h+O(h^2),
\label{eq:R22_local_product}
\end{equation}
and the ELIP approximation becomes
\begin{equation}
\frac{P_{\mathrm{aprox}}}{a+b}
=
\pi+
\left[(4-\pi)-R_{2,2}(1)\right]h
+O(h^2).
\label{eq:R22_local_ELIP}
\end{equation}
Matching the coefficient of $h$ in Equation~\eqref{eq:R22_local_ELIP} with the exact coefficient
$\pi/4$ in Equation~\eqref{eq:exact_circle_series_RL} gives
\begin{equation}
(4-\pi)-R_{2,2}(1)=\frac{\pi}{4},
\label{eq:R22_match_step}
\end{equation}
and therefore
\begin{equation}
\boxed{
R_{2,2}(1)=4-\frac{5\pi}{4}.
}
\label{eq:R22_T1}
\end{equation}
For the particular $R_{2,2}$ structure, this condition is explicitly
\begin{equation}
\frac{\pi-3+a_1+a_2}
{1+b_1+b_2}
=
4-\frac{5\pi}{4}.
\label{eq:R22_T1_coefficients}
\end{equation}

The notation $[T_1/F_3]$ is intended to make the construction transparent. The rational core
contains four unknown coefficients, $a_1,a_2,b_1,b_2$. One algebraic degree of freedom is
removed by the exact Taylor condition in Equation~\eqref{eq:R22_T1_coefficients}, while the
remaining three coefficients are intentionally left free and are determined by minimax fitting over
the complete interval $0<h<1$. This is a deliberate compromise: some of the available degrees
of freedom are used to enforce exact local behavior near the circle, while the others are retained
to control the global maximum error. Fully Taylor-determined versions, in which no free parameter
is retained, are considered separately in the following subsection.

For example, Equation~\eqref{eq:R22_T1_coefficients} may be used to eliminate $a_2$:
\begin{equation}
a_2=
\left(4-\frac{5\pi}{4}\right)(1+b_1+b_2)
-(\pi-3)-a_1.
\label{eq:R22_eliminate_a2}
\end{equation}
Thus only $a_1$, $b_1$, and $b_2$ need to be adjusted numerically. A minimax fit gives
\begin{equation}
\begin{aligned}
a_1&=1.67736246, &
a_2&=0.431876259,\\
b_1&=14.97383468, &
b_2&=14.85559359.
\end{aligned}
\label{eq:R22_fit_parameters_revised}
\end{equation}
The rational function is therefore
\begin{equation}
R_{2,2}(z)=
\frac{\pi-3+1.67736246z+0.431876259z^2}
{1+14.97383468z+14.85559359z^2},
\label{eq:R22_fit}
\end{equation}
and, because this member is short enough to be written explicitly, its complete perimeter formula is
\begin{equation}
\boxed{
P_{\mathrm{aprox}}
=
(a+b)
\left[
\pi+(4-\pi)h
-hz
\frac{\pi-3+1.67736246z+0.431876259z^2}
{1+14.97383468z+14.85559359z^2}
\right],
\qquad z=1-h.
}
\label{eq:R22_complete_revised}
\end{equation}
Its maximum relative error is approximately
\begin{equation}
E_{\max}\simeq3.3981\ \mathrm{ppm}.
\label{eq:R22_error_revised}
\end{equation}

A simpler rationalized implementation is
\begin{equation}
\boxed{
P_{\mathrm{aprox}}
=
(a+b)
\left[
\pi+(4-\pi)h
-hz
\frac{\pi-3+\frac53z+\frac37z^2}
{1+\frac{134}{9}z+\frac{59}{4}z^2}
\right],
\qquad z=1-h,
}
\label{eq:R22_rat}
\end{equation}
whose maximum relative error is approximately $3.5680\ \mathrm{ppm}$.

Adding the first exact logarithmic term in this family gives

\begin{equation}
P_{22,1\log}
=
(a+b)
\left[
\pi+(4-\pi)h
-hzR_{2,2}(z)
+\frac18z^2\ln\frac1z
\right].
\label{eq:R22_1log}
\end{equation}

The corresponding first-order Taylor condition changes because the logarithmic block itself
contributes to the coefficient of $h$. Indeed,
\begin{equation}
z^2\ln\frac1z
=
(1-h)^2\left[-\ln(1-h)\right]
=
h-\frac32h^2+O(h^3),
\label{eq:L1_circle_expansion}
\end{equation}
so the coefficient of $h$ in Equation~\eqref{eq:R22_1log}, after normalization by $a+b$, is
\begin{equation}
(4-\pi)-R_{2,2}(1)+\frac18.
\label{eq:R22_1log_hcoef}
\end{equation}
Matching this coefficient to the exact value $\pi/4$ gives
\begin{equation}
\boxed{
R_{2,2}(1)=\frac{33}{8}-\frac{5\pi}{4}.
}
\label{eq:R22_1log_T1}
\end{equation}
Again, one exact Taylor condition is imposed and three rational coefficients remain available for
global minimax fitting, hence the notation
\begin{equation}
\ELIP\;R_{2,2}+1\log[T_1/F_3].
\label{eq:R22_1log_name_revised}
\end{equation}
The resulting three-free-parameter fit reaches approximately $0.4841\ \mathrm{ppm}$. A particularly
compact integer-rational implementation is
\begin{equation}
\boxed{
P_{22,1\log}^{(67)}
=
(a+b)
\left[
\pi+(4-\pi)h
-hz\frac{67(\pi-3)+40z+22z^2}
{67+234z+60z^2}
+\frac18z^2\ln\frac1z
\right].
}
\label{eq:R22_1log_rat}
\end{equation}
which gives approximately $0.9041\ \mathrm{ppm}$.

The same design principle is used in the subsequent fitted members: enough Taylor conditions are
imposed to preserve exact local information, but three or four degrees of freedom are deliberately
retained for global minimax optimization. In the present notation, the progression is
\begin{equation}
\ELIP\;R_{2,2}+2\log[T_1/F_3],
\label{eq:name_R22L2}
\end{equation}
\begin{equation}
\ELIP\;R_{3,2}+2\log[T_2/F_3],
\label{eq:name_R32L2}
\end{equation}
\begin{equation}
\ELIP\;R_{4,2}+3\log[T_2/F_4],
\label{eq:name_R42L3}
\end{equation}
which yield approximately $0.2707$, $0.04395$, and $0.003236\ \mathrm{ppm}$, respectively.
Thus, the notation records explicitly how the available coefficients are divided between exact
Taylor information and numerical freedom.

The most accurate member of this partially constrained fitted sequence is
\begin{equation}
R_{5,2}(z)=
\frac{\pi-3+a_1z+a_2z^2+a_3z^3+a_4z^4+a_5z^5}
{1+b_1z+b_2z^2}.
\label{eq:R52_structure}
\end{equation}
This rational core has seven adjustable coefficients. Four of them are fixed through exact Taylor
conditions,
\begin{equation}
R(1)=\frac{4395}{1024}-\frac{5\pi}{4},
\label{eq:R52_T1}
\end{equation}
\begin{equation}
R'(1)=-\frac{7401}{2048}+\frac{81\pi}{64},
\label{eq:R52_T2}
\end{equation}
\begin{equation}
R''(1)=\frac{25465}{3072}-\frac{325\pi}{128},
\label{eq:R52_T3}
\end{equation}
\begin{equation}
R'''(1)=-\frac{48895}{2048}+\frac{62475\pi}{8192},
\label{eq:R52_T4}
\end{equation}
while three degrees of freedom are intentionally retained for the global minimax fit. The resulting
formula is therefore denoted
\begin{equation}
\ELIP\;R_{5,2}+3\log[T_4/F_3].
\label{eq:R52_name_revised}
\end{equation}
Its coefficients are
\begin{equation}
\begin{aligned}
a_1&=0.434456960000000, &
a_2&=0.425364514273400, &
a_3&=0.318693145897231,\\
a_4&=0.128019235361406, &
a_5&=0.00374776439801632,\\
b_1&=2.309109350000000, &
b_2&=0.668613450000000,
\end{aligned}
\label{eq:R52_cert_params}
\end{equation}
and its complete form is
\begin{equation}
P_{52}
=
(a+b)
\left[
\pi+(4-\pi)h
-hzR_{5,2}(z)
+L_3(z)
\right].
\label{eq:R52_cert_formula}
\end{equation}
The rigorous error certificate for this formula is discussed separately in
Section~\ref{sec:certification}.

\subsection{Fully Taylor-determined rational cores}\label{sec:zero_free}
A natural experiment is to remove numerical fitting from the rational core entirely. For an $R_{n,m}$ architecture there are $n+m$ coefficients after the constant term $\pi-3$ has been fixed. If all of them are determined by consecutive Maclaurin conditions, the notation becomes $T_{n+m}/F_0$.

The calculation is linear when organized at the level of series coefficients. Define the target rational residual by rearranging Equation~\eqref{eq:RL_general}:
\begin{equation}
R_{\mathrm{target}}(h)=\frac{\pi+(4-\pi)h+L_q(1-h)-G_{\mathrm{exact}}(h)}{h(1-h)}.
\label{eq:Rtarget}
\end{equation}
Expand it around $h=0$ as
\begin{equation}
R_{\mathrm{target}}(h)=r_0+r_1h+r_2h^2+\cdots+r_{n+m-1}h^{n+m-1}+O(h^{n+m}).
\label{eq:Rtarget_series}
\end{equation}
Next rewrite the denominator and numerator of Equation~\eqref{eq:Rnm_general} in powers of $h=1-z$:
\begin{equation}
D(h)=1+b_1(1-h)+\cdots+b_m(1-h)^m,
\label{eq:D_h}
\end{equation}
\begin{equation}
N(h)=\pi-3+a_1(1-h)+\cdots+a_n(1-h)^n.
\label{eq:N_h}
\end{equation}
Exact matching through the required order is equivalent to
\begin{equation}
N(h)-D(h)R_{\mathrm{target}}(h)=O(h^{n+m}).
\label{eq:linear_match}
\end{equation}
Because the $r_k$ in Equation~\eqref{eq:Rtarget_series} are already known exactly from $\pi$ and rational numbers, coefficient comparison in Equation~\eqref{eq:linear_match} gives a linear system for $a_i,b_j$. Solving this system in exact symbolic arithmetic yields coefficients in $\mathbb{Q}(\pi)$; no nonlinear fit is required.

\paragraph{Worked example: the fully Taylor-determined
$\ELIP\;R_{2,2}+1\log[T_4/F_0]$ formula.}
To make the construction above more explicit, we now develop one complete
example.  The purpose is to show, step by step, how all four coefficients of
the rational core can be obtained algebraically from the Maclaurin expansion,
without any numerical fitting.

For this case the normalized approximation is
\begin{equation}
\frac{P_{\mathrm{aprox}}}{a+b}
=
\pi+(4-\pi)h
-h(1-h)R_{2,2}(1-h)
+\frac18(1-h)^2
\ln\!\left(\frac{1}{1-h}\right),
\label{eq:R22L1_zero_start}
\end{equation}
with
\begin{equation}
R_{2,2}(z)
=
\frac{\pi-3+a_1z+a_2z^2}
{1+b_1z+b_2z^2},
\qquad z=1-h.
\label{eq:R22L1_zero_R}
\end{equation}
The four unknown quantities are therefore
$a_1$, $a_2$, $b_1$, and $b_2$.  Since the notation is
$[T_4/F_0]$, all four of them will be fixed by exact Taylor conditions,
leaving no free parameter for numerical optimization.

The exact normalized ellipse perimeter has the Maclaurin expansion
\begin{equation}
G_{\mathrm{exact}}(h)
=
\frac{P_{\mathrm{exact}}}{a+b}
=
\pi
+\frac{\pi}{4}h
+\frac{\pi}{64}h^2
+\frac{\pi}{256}h^3
+\frac{25\pi}{16384}h^4
+O(h^5).
\label{eq:R22L1_zero_exact_series}
\end{equation}
The logarithmic contribution in
Equation~\eqref{eq:R22L1_zero_start} is expanded independently:
\begin{equation}
\frac18(1-h)^2
\ln\!\left(\frac{1}{1-h}\right)
=
\frac18h
-\frac{3}{16}h^2
+\frac{1}{24}h^3
+\frac{1}{96}h^4
+O(h^5).
\label{eq:R22L1_zero_log_series}
\end{equation}

Next, instead of expanding the four rational coefficients directly, it is
convenient to first determine the Taylor series that the rational function
itself must reproduce.  Write
\begin{equation}
R_{2,2}(1-h)
=
r_0+r_1h+r_2h^2+r_3h^3+O(h^4).
\label{eq:R22L1_zero_rseries}
\end{equation}
Since
\begin{equation}
-h(1-h)R_{2,2}(1-h)
=
-r_0h
+(r_0-r_1)h^2
+(r_1-r_2)h^3
+(r_2-r_3)h^4
+O(h^5),
\label{eq:R22L1_zero_product_series}
\end{equation}
matching Equations~\eqref{eq:R22L1_zero_start}--\eqref{eq:R22L1_zero_product_series}
with the exact series in Equation~\eqref{eq:R22L1_zero_exact_series}
gives the four successive equations
\begin{equation}
(4-\pi)+\frac18-r_0=\frac{\pi}{4},
\label{eq:R22L1_zero_match1}
\end{equation}
\begin{equation}
-\frac{3}{16}+r_0-r_1=\frac{\pi}{64},
\label{eq:R22L1_zero_match2}
\end{equation}
\begin{equation}
\frac{1}{24}+r_1-r_2=\frac{\pi}{256},
\label{eq:R22L1_zero_match3}
\end{equation}
\begin{equation}
\frac{1}{96}+r_2-r_3=\frac{25\pi}{16384}.
\label{eq:R22L1_zero_match4}
\end{equation}
Solving them sequentially gives
\begin{equation}
r_0=
\frac{33}{8}-\frac{5\pi}{4},
\label{eq:R22L1_zero_r0}
\end{equation}
\begin{equation}
r_1=
\frac{63}{16}-\frac{81\pi}{64},
\label{eq:R22L1_zero_r1}
\end{equation}
\begin{equation}
r_2=
\frac{191}{48}-\frac{325\pi}{256},
\label{eq:R22L1_zero_r2}
\end{equation}
\begin{equation}
r_3=
\frac{383}{96}-\frac{20825\pi}{16384}.
\label{eq:R22L1_zero_r3}
\end{equation}
Thus, the Taylor conditions first determine exactly the local series that
$R_{2,2}(z)$ must have at the circular endpoint $z=1$.

It remains to convert these four exact series coefficients into the four
coefficients of the rational function.  Rewriting its numerator and
denominator in powers of $h$ gives
\begin{equation}
D(h)
=
1+b_1(1-h)+b_2(1-h)^2
=
d_0+d_1h+d_2h^2,
\label{eq:R22L1_zero_D}
\end{equation}
where
\begin{equation}
d_0=1+b_1+b_2,
\qquad
d_1=-b_1-2b_2,
\qquad
d_2=b_2,
\label{eq:R22L1_zero_dcoeff}
\end{equation}
and
\begin{equation}
N(h)
=
\pi-3+a_1(1-h)+a_2(1-h)^2
=
n_0+n_1h+n_2h^2,
\label{eq:R22L1_zero_N}
\end{equation}
with
\begin{equation}
n_0=\pi-3+a_1+a_2,
\qquad
n_1=-a_1-2a_2,
\qquad
n_2=a_2.
\label{eq:R22L1_zero_ncoeff}
\end{equation}

Because
\begin{equation}
R_{2,2}(1-h)=\frac{N(h)}{D(h)},
\label{eq:R22L1_zero_ratio}
\end{equation}
the desired matching condition is equivalent to
\begin{equation}
N(h)
=
D(h)
\left(
r_0+r_1h+r_2h^2+r_3h^3
\right)
+O(h^4).
\label{eq:R22L1_zero_identity}
\end{equation}
Comparison of equal powers of $h$ yields
\begin{equation}
n_0=d_0r_0,
\label{eq:R22L1_zero_linear0}
\end{equation}
\begin{equation}
n_1=d_0r_1+d_1r_0,
\label{eq:R22L1_zero_linear1}
\end{equation}
\begin{equation}
n_2=d_0r_2+d_1r_1+d_2r_0,
\label{eq:R22L1_zero_linear2}
\end{equation}
\begin{equation}
0=d_0r_3+d_1r_2+d_2r_1.
\label{eq:R22L1_zero_linear3}
\end{equation}
Equations~\eqref{eq:R22L1_zero_linear0}--\eqref{eq:R22L1_zero_linear3}
form a linear algebraic system for
$a_1,a_2,b_1,b_2$.  This is the key point: once the four exact
$r_k$ values have been obtained from the Maclaurin expansion, no nonlinear
optimization is needed.

For completeness, the exact solution may be written compactly by defining
\begin{equation}
\Delta
=
28575\pi^2-464772\pi+1177600.
\label{eq:R22L1_zero_Delta}
\end{equation}
Then
\begin{equation}
a_1=
\frac{
-18962550\pi^3
+175652430\pi^2
-542002193\pi
+557088000
}{\Delta},
\label{eq:R22L1_zero_a1}
\end{equation}
\begin{equation}
a_2=
\frac{
141975\pi^3
+3050547\pi^2
-23217467\pi
+38429952
}{\Delta},
\label{eq:R22L1_zero_a2}
\end{equation}
\begin{equation}
b_1=
\frac{
-154350\pi^2
+2395968\pi
-6004736
}{\Delta},
\label{eq:R22L1_zero_b1}
\end{equation}
\begin{equation}
b_2=
\frac{
15159375\pi^2
-94842300\pi
+148338688
}{\Delta}.
\label{eq:R22L1_zero_b2}
\end{equation}
Numerically,
\begin{equation}
\begin{aligned}
a_1&\simeq0.366604358879268, &
a_2&\simeq0.126976933868143,\\
b_1&\simeq1.906590738158613, &
b_2&\simeq0.301209711663339.
\end{aligned}
\label{eq:R22L1_zero_decimal}
\end{equation}

Substitution of these coefficients into
Equation~\eqref{eq:R22L1_zero_start} gives
$\ELIP\;R_{2,2}+1\log[T_4/F_0]$.  By construction,
\begin{equation}
\frac{P_{\mathrm{aprox}}}{a+b}
-
G_{\mathrm{exact}}(h)
=
O(h^5),
\qquad h\rightarrow0,
\label{eq:R22L1_zero_order}
\end{equation}
so the approximation reproduces the exact perimeter series through
$h^4$.  All four rational coefficients are therefore consequences of the
Taylor conditions; none has been fitted numerically.  A dense numerical
evaluation over $0<h<1$ gives a maximum relative error of approximately
$7.0200$ ppm, as reported in Table~\ref{tab:zero_free}.

Table~\ref{tab:zero_free} summarizes the sequence obtained by applying this construction from the smallest rational core to $R_{5,3}+3\log$. These maxima are numerical scans and are included to show the trend, not as certificates.

\begin{table}[H]
\centering
\footnotesize
\renewcommand{\arraystretch}{1.16}
\caption{Fully Taylor-determined rational cores.}
\label{tab:zero_free}
\begin{tabular}{lccc}
\toprule
Formula & Taylor order & Free rational parameters & $E_{\max}$ (ppm)\\
\midrule
$\ELIP\;R_{2,2}[T_4/F_0]$ & $h^4$ & 0 & 41.9866\\
$\ELIP\;R_{2,2}+1\log[T_4/F_0]$ & $h^4$ & 0 & 7.0200\\
$\ELIP\;R_{2,2}+2\log[T_4/F_0]$ & $h^4$ & 0 & 4.2880\\
$\ELIP\;R_{3,2}+2\log[T_5/F_0]$ & $h^5$ & 0 & 1.0508\\
$\ELIP\;R_{4,2}+3\log[T_6/F_0]$ & $h^6$ & 0 & 0.22393\\
$\ELIP\;R_{5,2}+3\log[T_7/F_0]$ & $h^7$ & 0 & 0.075706\\
$\ELIP\;R_{5,3}+3\log[T_8/F_0]$ & $h^8$ & 0 & 0.025477\\
\bottomrule
\end{tabular}
\end{table}

For the two highest-degree cases whose coefficients are fully determined by exact Taylor matching (i.e., with no free parameters), the coefficients are given by:
\begin{equation}
\begin{aligned}
R_{5,2}:\quad &(a_1,a_2,a_3,a_4,a_5)=(0.3150934886423322,0.2815750673412721,\\
&0.2039979780492760,0.0614905268148402,0.0011348514691154),\\
&(b_1,b_2)=(1.4667445437665716,0.2863535185638652),
\end{aligned}
\label{eq:R52_zero_params}
\end{equation}
\begin{equation}
\begin{aligned}
R_{5,3}:\quad &(a_1,a_2,a_3,a_4,a_5)=(0.4197742928969318,0.4422532000280925,\\
&0.3361123798694394,0.1479444397542482,0.01479665860551144),\\
&(b_1,b_2,b_3)=(2.205678177300706,0.8637674197543956,0.04690504193927536).
\end{aligned}
\label{eq:R53_zero_params}
\end{equation}
The displayed decimals are implementation values; the exact symbolic solutions are rational functions of $\pi$. To make the fully algebraic character of the $R_{5,3}$ core explicit,
the eight coefficients reported numerically in Eq.~\eqref{eq:R53_zero_params} can be written
exactly as rational functions of $\pi$. Define

\begin{equation}
\begin{aligned}
\Delta ={}&
65772081577125\pi^{3}
-13588421557756500\pi^{2}
+110225586548152320\pi \\
&-214210269553885184 .
\end{aligned}
\label{eq:R53_exact_delta}
\end{equation}

Then the five numerator coefficients are

\begin{equation}
\begin{aligned}
a_1={}&
\frac{1}{768000\Delta}
\Big(
-2575158513235247352105000\pi^{4}
+250436150359939914305079375\pi^{3}
\\
&\qquad
-2206737723616652919855462300\pi^{2}
+6769583900607101907603448320\pi
\\
&\qquad
-7001895419515260595356041216
\Big),
\\[3pt]
a_2={}&
\frac{1}{384000\Delta}
\Big(
-970241513306663831310000\pi^{4}
+23322323610048894885664875\pi^{3}
\\
&\qquad
-175545198467187256180185300\pi^{2}
+532751531468312688460492800\pi
\\
&\qquad
-569754708747483109799231488
\Big),
\\[3pt]
a_3={}&
\frac{1}{192000\Delta}
\Big(
3783033694991622802500\pi^{4}
+353726032639417370186625\pi^{3}
\\
&\qquad
-8064205316346185927802000\pi^{2}
+39718113733264010784975360\pi
\\
&\qquad
-56523847107622410559422464
\Big).
\end{aligned}
\label{eq:R53_exact_a123}
\end{equation}

\begin{equation}
\begin{aligned}
a_4={}&
\frac{1}{128000\Delta}
\Big(
-262313325373210650000\pi^{4}
-19054842336504991382625\pi^{3}
\\
&\qquad
-3279906924829433809085700\pi^{2}
+21201364277442120499752960\pi
\\
&\qquad
-33618295011611685081317376
\Big),
\\[3pt]
a_5={}&
\frac{1}{256000\Delta}
\Big(
37427698213418205000\pi^{4}
+118787995715378766753375\pi^{3}
\\
&\qquad
-3038011793519873764740300\pi^{2}
+15565540376162329814684160\pi
\\
&\qquad
-22603431984164045986988032
\Big).
\end{aligned}
\label{eq:R53_exact_a45}
\end{equation}

The three denominator coefficients are

\begin{equation}
\begin{aligned}
b_1={}&
\frac{1}{25\Delta}
\Big(
-19724850981210375\pi^{3}
+3727048063637889900\pi^{2}
\\
&\qquad
-29136269906176849920\pi
+55361430198379610112
\Big),
\\[3pt]
b_2={}&
\frac{1}{25\Delta}
\Big(
65948701831397760375\pi^{3}
-6214003414361822979900\pi^{2}
\\
&\qquad
+37086443925480555217920\pi
-57225568019472829120512
\Big),
\\[3pt]
b_3={}&
\frac{1}{25\Delta}
\Big(
51306891451074133875\pi^{3}
-1166190775849699746300\pi^{2}
\\
&\qquad
+5807917804903896407040\pi
-8327105965908027244544
\Big).
\end{aligned}
\label{eq:R53_exact_b123}
\end{equation}

Thus, none of the eight coefficients of the rational core
$R_{5,3}$ is obtained by numerical fitting: all of them are fixed
algebraically by the eight exact Taylor-matching conditions. Their
decimal values in Eq.~\eqref{eq:R53_zero_params} are merely numerical representations of
the exact expressions above. The exact symbolic expressions are intended for analytical characterization. For numerical implementation, the decimal coefficients in Eq.~\eqref{eq:R53_zero_params} are recommended.

The later beta-corrected approximation
$\mathrm{ELIP}\;R_{5,3}+3\log[T_8/F_0;\beta_1]$ uses this same
fully Taylor-determined core. Its additional correction is

\begin{equation}
P_{53,\beta_1}
=
P_{53,0}
\left(
1-A_{\beta}h^{12}z^{7/4}
\right),
\qquad
A_{\beta}
=
4.8665743356537865\times10^{-6},
\end{equation}

where $A_{\beta}$ is the only coefficient determined by numerical
minimax fitting. This correction reduces the maximum absolute
relative error to approximately $0.0004512$ ppm.

\subsection{Beta-shaped residuals and compensation}\label{sec:beta}
\begin{figure}[H]
\centering
\includegraphics[width=0.78\linewidth]{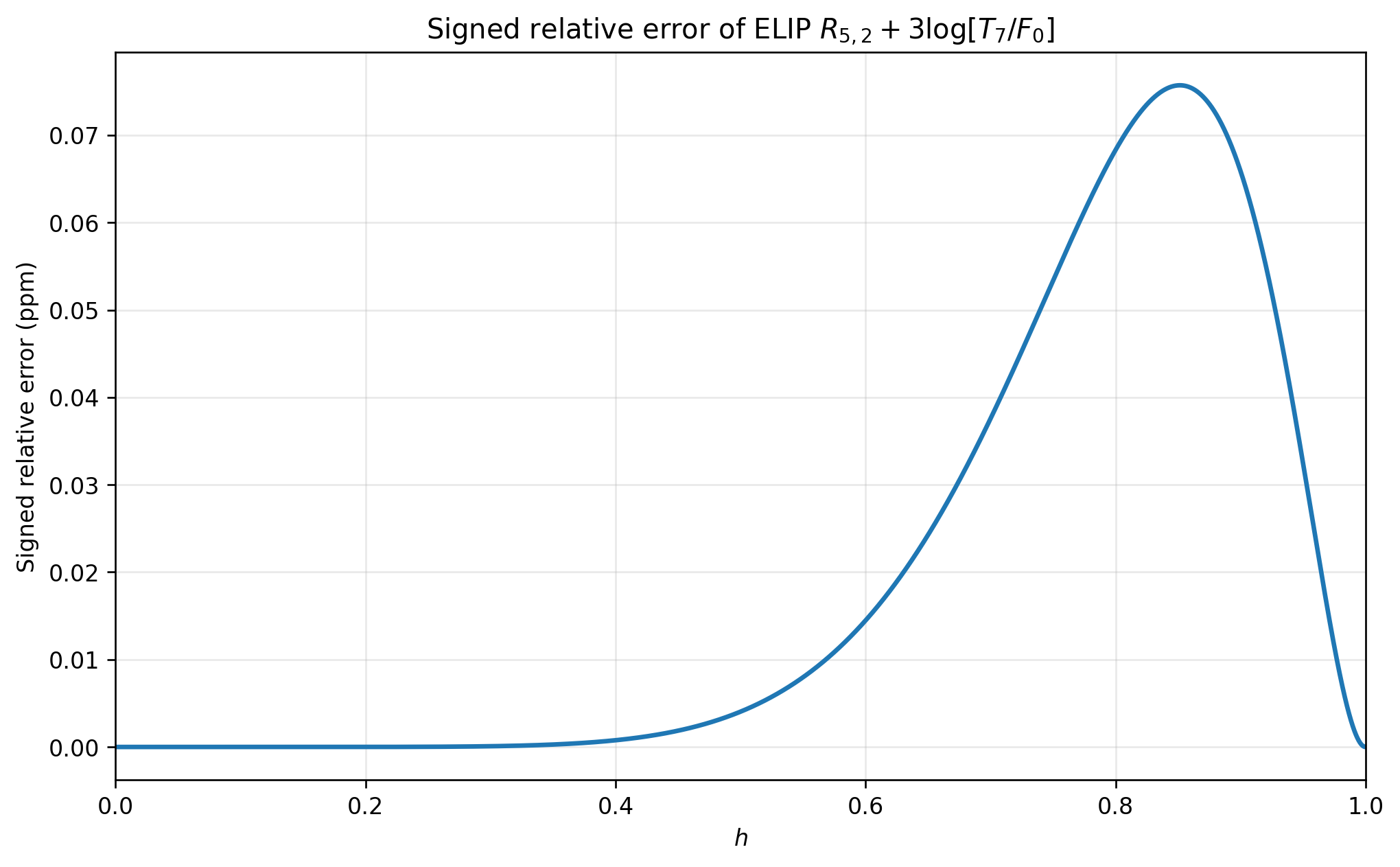}
\caption{Signed relative error of the fully Taylor-determined $\ELIP\;R_{5,2}+3\log[T_7/F_0]$ core.}
\label{fig:beta52}
\end{figure}

\begin{figure}[H]
\centering
\includegraphics[width=0.78\linewidth]{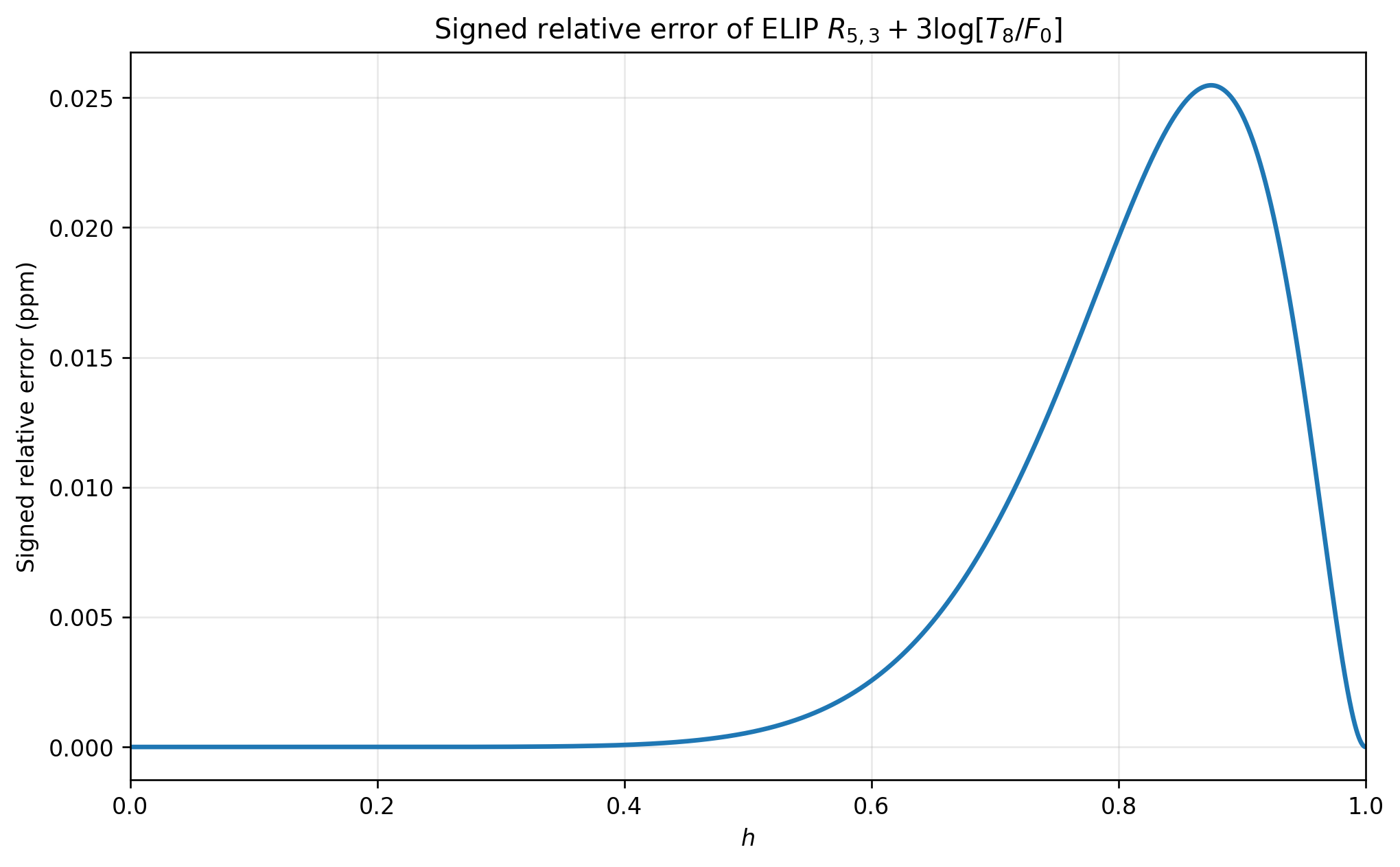}
\caption{Signed relative error of the fully Taylor-determined $\ELIP\;R_{5,3}+3\log[T_8/F_0]$ core.}
\label{fig:beta53}
\end{figure}

Once every rational coefficient is consumed by Taylor restrictions, no minimax freedom remains to flatten the residual away from the circle. For $R_{5,2}$ and $R_{5,3}$, the signed errors are nevertheless remarkably smooth and resemble a single beta-distribution lobe, as shown in Figures~\ref{fig:beta52} and \ref{fig:beta53}.

This shape suggests a compensator that vanishes at both endpoints and begins only beyond the exact Maclaurin order. For the $R_{5,2}$ core we use
\begin{equation}
P_{52,\beta2}=P_{52,0}\left[1-A_{52}\frac{h^8z^2}{1+B_{52}z}\right],
\label{eq:R52_beta2}
\end{equation}
with
\begin{equation}
A_{52}=3.154299880686168\times10^{-5},\qquad B_{52}=10.73826918095434.
\label{eq:R52_beta2_params}
\end{equation}
The resulting notation is $\ELIP\;R_{5,2}+3\log[T_7/F_0;\beta_2]$, and the numerical maximum is approximately
\begin{equation}
E_{\max}\simeq0.0016331\ppm.
\label{eq:R52_beta2_error}
\end{equation}
Near the circle, the compensator is $O(h^8)$, so the high-order local accuracy of the fully constrained core is retained to the intended order.

For $R_{5,3}$, the analogous two-parameter correction is
\begin{equation}
P_{53,\beta2}=P_{53,0}\left[1-A_{53}\frac{h^9z^2}{1+B_{53}z}\right],
\label{eq:R53_beta2}
\end{equation}
where
\begin{equation}
A_{53}=1.867101678559634\times10^{-5},\qquad B_{53}=20.62831105372712.
\label{eq:R53_beta2_params}
\end{equation}
This formula is referred to as $\ELIP\;R_{5,3}+3\log[T_8/F_0;\beta_2]$ and its numerical minimax error is
\begin{equation}
E_{\max}\simeq0.00096813\ppm.
\label{eq:R53_beta2_error}
\end{equation}
The formerly single positive lobe is redistributed into three extrema of nearly equal magnitude and alternating sign, which is the expected minimax mechanism.

An even simpler compensator is obtained by fixing both exponents at values selected through numerical exploration and optimizing only the amplitude:
\begin{equation}
P_{53,\beta1}=P_{53,0}\left[1-Ah^{12}z^{7/4}\right].
\label{eq:R53_beta1}
\end{equation}
The fitted value is
\begin{equation}
A=4.8665743356537865\times10^{-6}.
\label{eq:R53_beta1_A}
\end{equation}
The correction in Equation~\eqref{eq:R53_beta1} begins at order $h^{12}$, so it does not alter the exact $h^8$ Maclaurin matching of the base. Its numerical maximum relative error is about $0.00045119\ppm$, and Section~\ref{sec:certification} gives a rigorous global bound. The complete formula named as $\ELIP\;R_{5,3}+3\log[T_8/F_0;\beta_1]$ is

\begin{equation}
\boxed{
P_{53,\beta_1}
=
(a+b)
\left[
\pi+(4-\pi)h
-hzR_{5,3}(z)
+
\left(
\frac{z^2}{8}
+\frac{3z^3}{32}
+\frac{75z^4}{1024}
\right)
\ln\!\left(\frac{1}{z}\right)
\right]
\left[
1-Ah^{12}z^{7/4}
\right],
}
\label{eq:R53_beta1_explicit}
\end{equation}
where the fully Taylor-determined rational function
$R_{5,3}(z)$ and its exact algebraic coefficients are calculated in Equations~\eqref{eq:R53_exact_a123}-\eqref{eq:R53_exact_b123}, and its numerical approximation values are given in
Equation~\eqref{eq:R53_zero_params}. The only numerically fitted
parameter in Equation~\eqref{eq:R53_beta1} is the correction
amplitude $A$.

\subsection{A purely rational low-complexity variant}\label{sec:R33}
One additional formula is useful because it contains no logarithm, exponential, or external beta factor. Let
\begin{equation}
R_{3,3}(z)=\frac{\pi-3+a_1z+a_2z^2+a_3z^3}{1+b_1z+b_2z^2+b_3z^3}.
\label{eq:R33_structure}
\end{equation}
Five algebraic conditions are imposed: four exact Maclaurin matches through $h^4$ and a fifth condition making the leading $h^5$ relative-error coefficient equal to Cantrell's. Only one parameter remains free. The minimax value is close to $b_3=28.98375$, and the particularly simple choice
\begin{equation}
b_3=29
\label{eq:R33_b3}
\end{equation}
leads to
\begin{equation}
\begin{aligned}
a_1&=3.671671646844034, & a_2&=6.166434371309410, & a_3&=0.371569742207349,\\
b_1&=29.624958210437996, & b_2&=82.155404023470153, & b_3&=29.
\end{aligned}
\label{eq:R33_params}
\end{equation}
The resulting $\ELIP\;R_{3,3}[T_5/F_1]$ has a maximum error near $1.59569\ppm$ and shares Cantrell's leading circular-limit relative error
\begin{equation}
\varepsilon(h)=-\frac{3}{131072}h^5+O(h^6).
\label{eq:R33_circle_error}
\end{equation}
This formula is included among the simple ELIP candidates because five of its six coefficients are algebraically determined and the remaining degree of freedom has been replaced by the nearby integer $29$.

\section{Numerical Evaluations and Discussion}\label{sec:numerics}

\subsection{Computational setup}
Exploratory minimax fits were performed with high-precision reference evaluations and dense grids enriched near both endpoints. Global optimization was followed by local refinement, and the signed-error extrema were independently localized. For the tables and figures in this section we evaluate every formula directly in the common variable $h$ and report absolute relative error in ppm according to Equation~\eqref{eq:error_definition}. The grid density is sufficient for plotting and comparative numerical estimates, but only the interval-certified values in Section~\ref{sec:certification} are treated as rigorous global bounds.

\subsection{Four compact ELIP formulas}\label{sec:simple_formulas}
Table~\ref{tab:simple} highlights four formulas intended to emphasize simplicity rather than the smallest possible error. Their architectures differ substantially, yet each uses only a small number of compact operations.

\begin{table}[H]
\centering
\footnotesize
\renewcommand{\arraystretch}{1.18}
\caption{Representative low-complexity ELIP formulas. The ``fitted'' column counts genuine numerical degrees of freedom in the parent construction; rationalized implementations replace some decimals by small fractions or integers.}
\label{tab:simple}
\begin{tabular}{p{5.7cm}ccc}
\toprule
Formula & Fitted parameters & Main operations & $E_{\max}$ (ppm)\\
\midrule
$\ELIP\;LPD[T_0/F_3]_{\mathrm{rat}}$ & 3 (during original fitting) & log, powers & 6.8404\\
$\ELIP\;R_{2,2}[T_1/F_3]_{\mathrm{rat}}$ & 3 (during original fitting) & arithmetic & 3.5680\\
$\ELIP\;R_{2,2}+1\log[T_1/F_3]_{\mathrm{rat}}$ & 3 (during original fitting) & rational + log & 0.9041\\
$\ELIP\;R_{3,3}[T_5/F_1]$, $b_3=29$ & 1 & arithmetic & 1.5957\\
\bottomrule
\end{tabular}
\end{table}

\begin{figure}[H]
\centering
\includegraphics[width=0.94\linewidth]{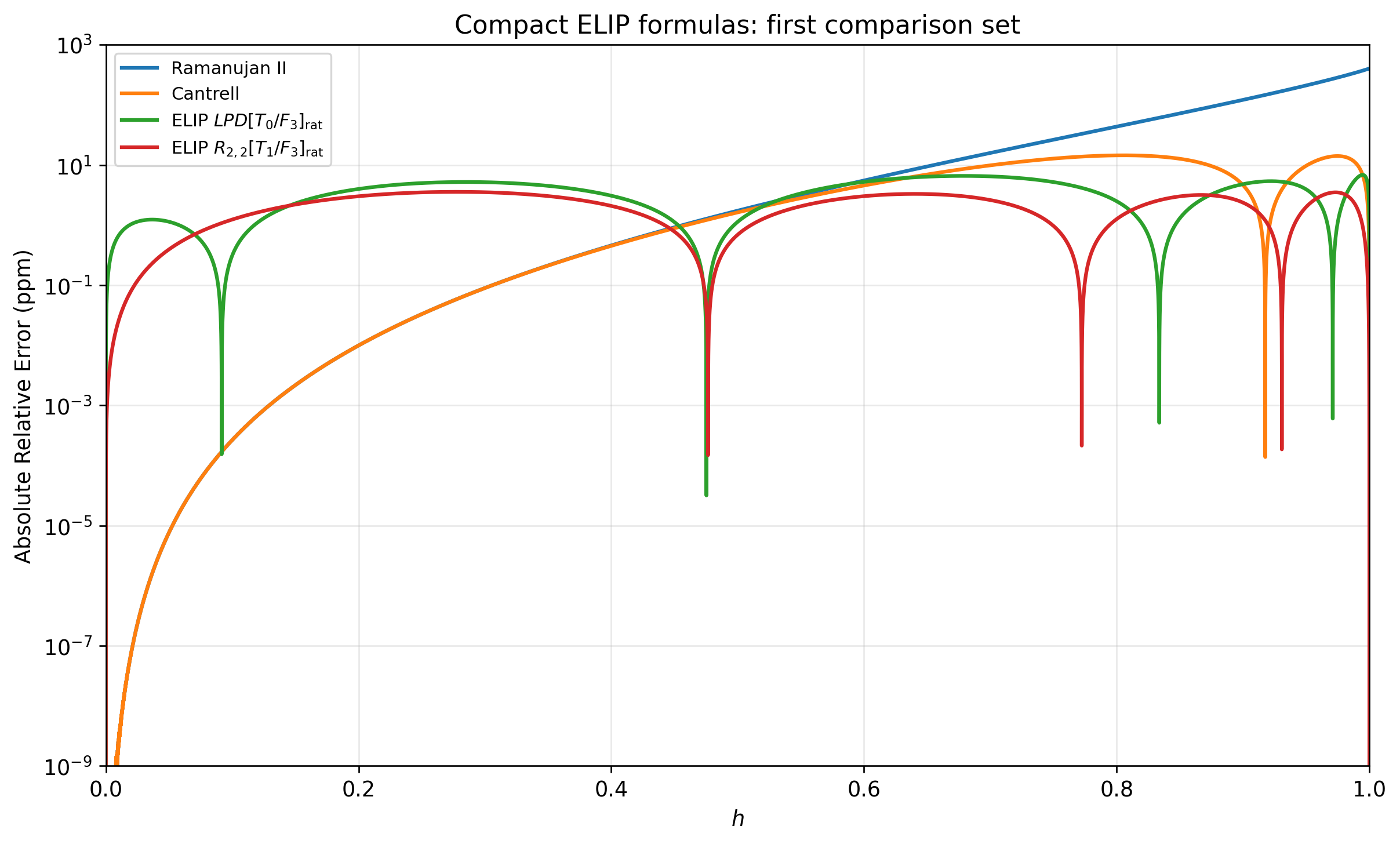}
\caption{Absolute relative error of Ramanujan II, Cantrell, $\ELIP\;LPD[T_0/F_3]_{\mathrm{rat}}$, and $\ELIP\;R_{2,2}[T_1/F_3]_{\mathrm{rat}}$. The vertical axis is shown from $10^{-9}$ to $10^3$ ppm.}
\label{fig:simplecomparisonA}
\end{figure}

\begin{figure}[H]
\centering
\includegraphics[width=0.94\linewidth]{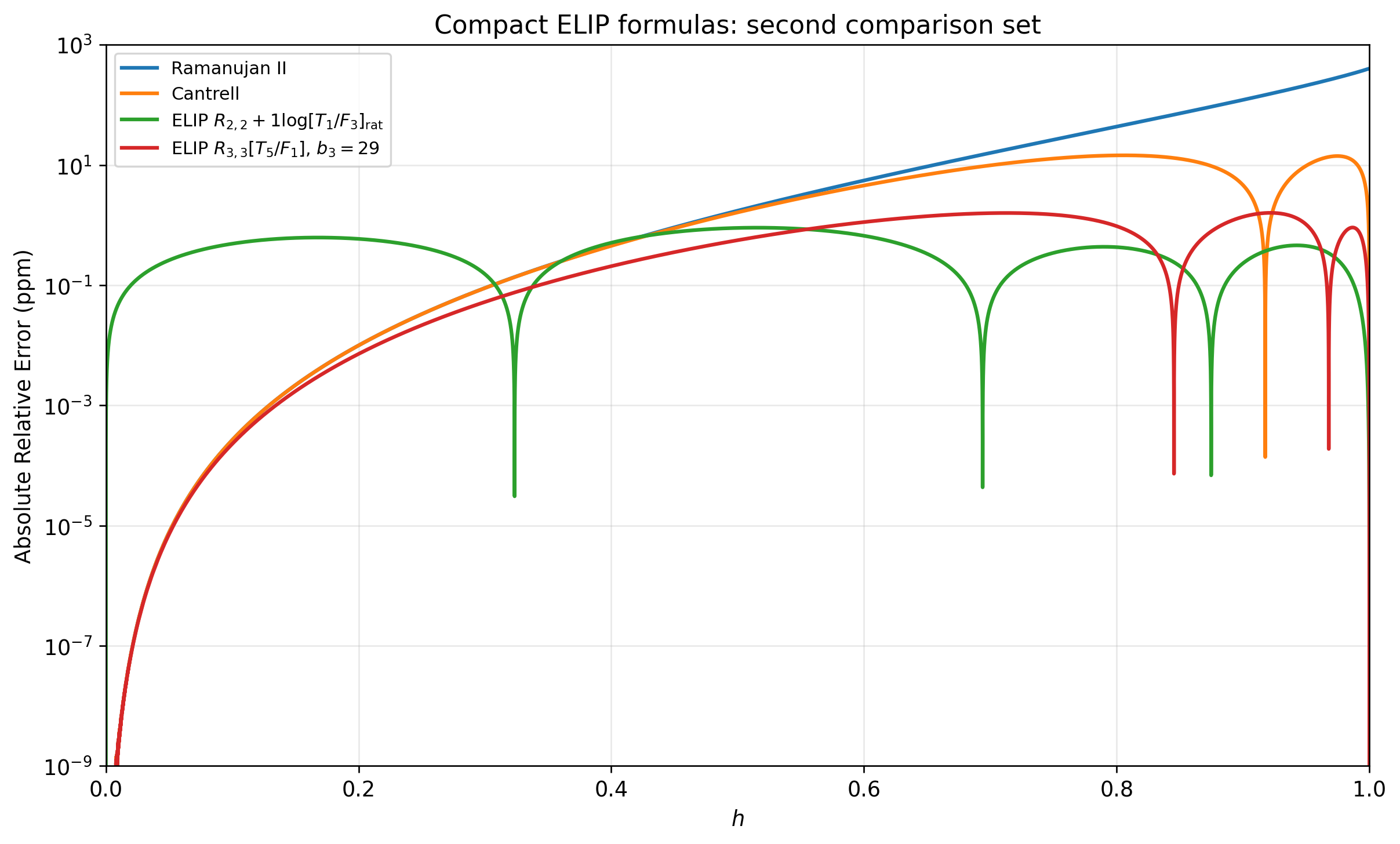}
\caption{Absolute relative error of Ramanujan II, Cantrell, $\ELIP\;R_{2,2}+1\log[T_1/F_3]_{\mathrm{rat}}$, and $\ELIP\;R_{3,3}[T_5/F_1]$ with $b_3=29$. The vertical axis is shown from $10^{-9}$ to $10^3$ ppm.}
\label{fig:simplecomparisonB}
\end{figure}

Figures~\ref{fig:simplecomparisonA} and~\ref{fig:simplecomparisonB} split the compact-formula comparison into two panels so that the curves are easier to distinguish. The rationalized $R_{2,2}+1\log$ member is particularly notable because it is sub-ppm with a very small rational function and one logarithm. The $R_{3,3}$ formula is somewhat less accurate globally, but it is purely rational and has the same leading circular-limit error coefficient as Cantrell.

\subsection{Three high-accuracy ELIP formulas}\label{sec:high_formulas}
The three principal high-accuracy formulas are:
\begin{equation}
\mathcal{E}_1=\ELIP\;R_{5,2}+3\log[T_4/F_3],
\label{eq:E1_name}
\end{equation}
\begin{equation}
\mathcal{E}_2=\ELIP\;R_{5,2}+3\log[T_7/F_0;\beta_2],
\label{eq:E2_name}
\end{equation}
\begin{equation}
\mathcal{E}_3=\ELIP\;R_{5,3}+3\log[T_8/F_0;\beta_1].
\label{eq:E3_name}
\end{equation}
The first achieves its accuracy through three fitted degrees of freedom in the rational core. The second determines the complete rational core by exact Taylor constraints and introduces only two fitted beta parameters. The third goes further: eight rational coefficients are Taylor-determined and only one external amplitude is fitted, while the correction exponents are fixed at $12$ and $7/4$.

Table~\ref{tab:highcomparison} compares these formulas with important reference methods over the full interval and over $0<h<0.3$. The near-circle column is especially informative because two formulas with similar global maxima may have dramatically different asymptotic quality.

\begin{table}[H]
\centering
\scriptsize
\renewcommand{\arraystretch}{1.18}
\caption{Numerical maximum absolute relative error over the full interval and near the circle. Certified inequalities are stated separately in Section~\ref{sec:certification}.}
\label{tab:highcomparison}
\begin{tabular}{p{5.4cm}cc}
\toprule
Method & $0<h<1$ (ppm) & $0<h<0.3$ (ppm)\\
\midrule
Ramanujan II & 402.3375 & 0.089682\\
Cantrell & 14.4608 & 0.089429\\
R2/2exp \cite{Ayala2026} & 0.5736 & 0.077831\\
Popenkov \cite{Popenkov2026} & 0.0158513 & $4.78\times10^{-5}$\\
UNDERWOOD F2, displayed 8-constant set & 0.01144 & 0.01143\\
UNDERWOOD U10, displayed 10-constant set & 0.0006512 & 0.0006400\\
$\mathcal{E}_1$: $R_{5,2}+3\log[T_4/F_3]$ & 0.0016226 & 0.0010483\\
$\mathcal{E}_2$: $R_{5,2}+3\log[T_7/F_0;\beta_2]$ & 0.0016331 & $3.48\times10^{-5}$\\
$\mathcal{E}_3$: $R_{5,3}+3\log[T_8/F_0;\beta_1]$ & 0.0004512 & $4.59\times10^{-6}$\\
\bottomrule
\end{tabular}
\end{table}

The global maxima of $\mathcal{E}_1$ and $\mathcal{E}_2$ are almost identical, but their circular-limit behavior is not. Formula $\mathcal{E}_2$ retains the exact $h^7$ matching of its zero-free rational core, so it is roughly two orders of magnitude more accurate than $\mathcal{E}_1$ on $0<h<0.3$. Formula $\mathcal{E}_3$ improves both criteria simultaneously: its full-range maximum is below the U10 value obtained from the ten-constant coefficient set used here, while its near-circle maximum is smaller by roughly two orders of magnitude.

\begin{figure}[H]
\centering
\includegraphics[width=0.94\linewidth]{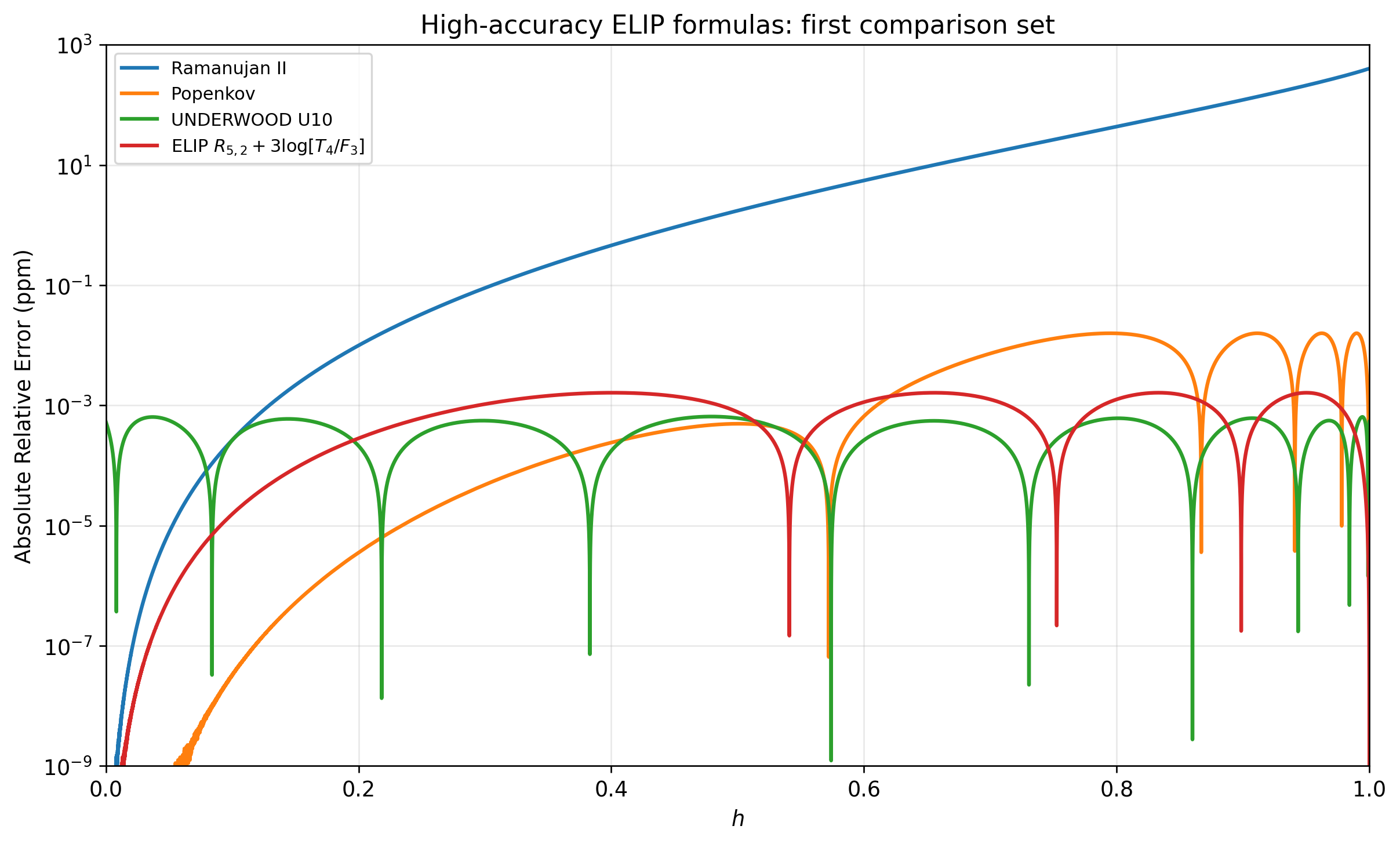}
\caption{Absolute relative error of Ramanujan II, Popenkov, UNDERWOOD U10, and $\mathcal{E}_1=\ELIP\;R_{5,2}+3\log[T_4/F_3]$. The vertical axis is shown from $10^{-9}$ to $10^3$ ppm.}
\label{fig:highcomparisonA}
\end{figure}

\begin{figure}[H]
\centering
\includegraphics[width=0.94\linewidth]{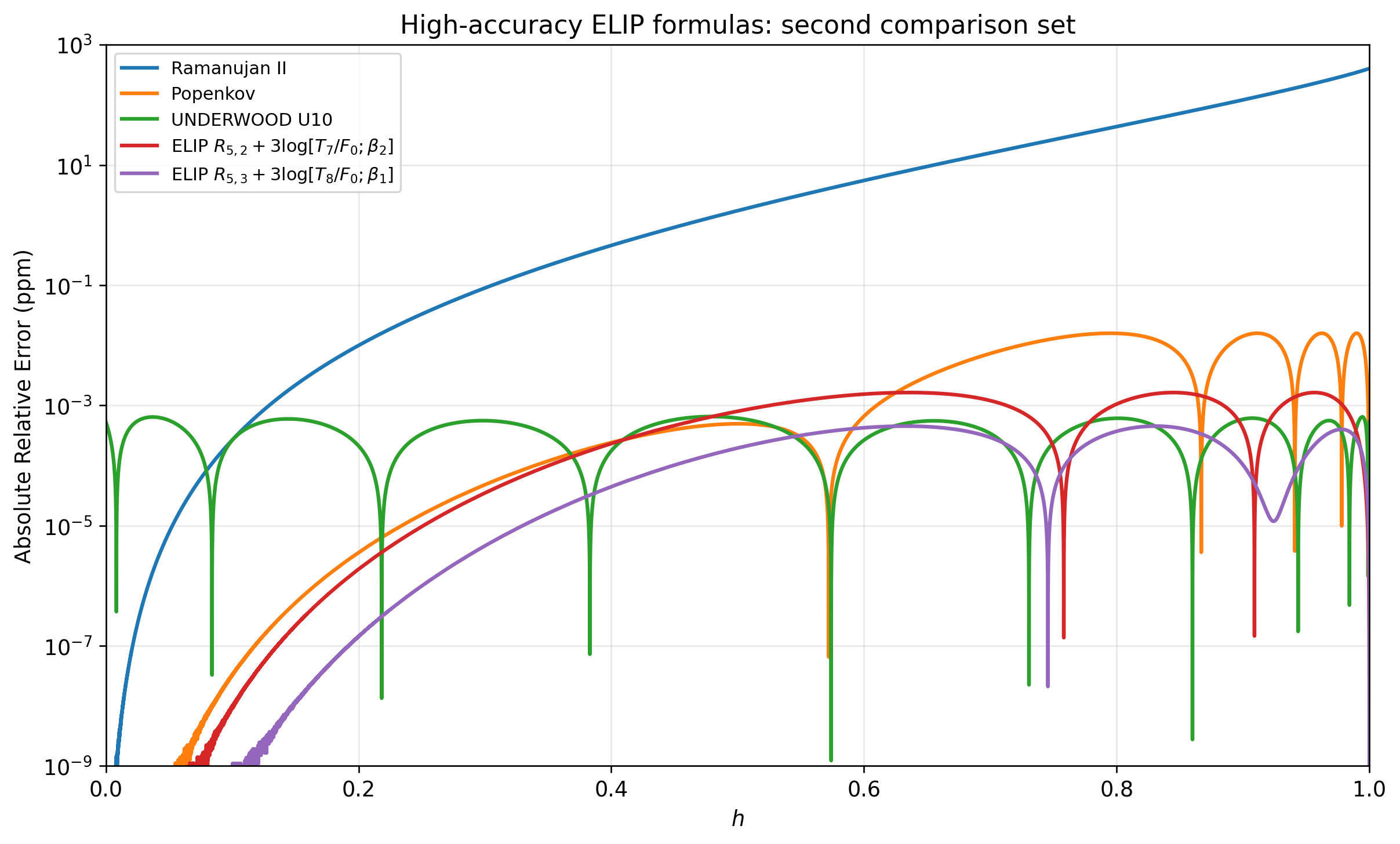}
\caption{Absolute relative error of Ramanujan II, Popenkov, UNDERWOOD U10, $\mathcal{E}_2=\ELIP\;R_{5,2}+3\log[T_7/F_0;\beta_2]$, and $\mathcal{E}_3=\ELIP\;R_{5,3}+3\log[T_8/F_0;\beta_1]$. The vertical axis is shown from $10^{-9}$ to $10^3$ ppm.}
\label{fig:highcomparisonB}
\end{figure}

Figures~\ref{fig:highcomparisonA} and~\ref{fig:highcomparisonB} provide the requested logarithmic comparisons, now split into two panels to avoid visual crowding. The curve for $\mathcal{E}_3$ falls very rapidly toward the circle because its base is exact through $h^8$ and its beta correction starts only at $h^{12}$.

\subsection{Interpretation of accuracy versus structure}
The numerical results reveal two distinct design philosophies. A classical minimax formula spends its adjustable coefficients to distribute error oscillations as uniformly as possible. A fully Taylor-constrained ELIP core spends those same coefficients on exact local structure; its residual can therefore be larger globally but smoother. The beta correction is effective precisely because the residual is simple enough to be compensated by a low-dimensional shape.

The comparison between $\mathcal{E}_1$ and $\mathcal{E}_2$ is illustrative. Their maximum errors differ by less than one percent, yet $\mathcal{E}_2$ has no fitted coefficient inside the rational core. Its fitted dimension is moved entirely into two external beta parameters, while the rational coefficients are fixed from the exact series. Thus, if ``simplicity'' is interpreted as the number of empirically fitted degrees of freedom rather than the number of stored constants, the two formulas have very different meanings despite nearly identical global accuracy.

The $R_{5,3}$ beta-one formula provides the clearest example. Its rational-logarithmic core is fully analytical, the external correction has only one fitted amplitude, and the fixed exponents were selected as simple values. The resulting maximum is around $4.5\times10^{-4}\ppm$, with much stronger circular-limit behavior than formulas that distribute approximation error more uniformly over the full interval.

\section{Certified Error Bounds}\label{sec:certification}

\subsection{Certification principle}
Dense numerical grids are useful for discovering candidate formulas, but they cannot prove that a narrow unsampled spike does not exist. For the selected ELIP formulas we therefore use Arb midpoint-radius arithmetic through \texttt{python-flint}. Let
\begin{equation}
\varepsilon(z)=\frac{G_{\mathrm{aprox}}(z)}{G_{\mathrm{exact}}(z)}-1.
\label{eq:cert_error}
\end{equation}
On an interval $I=[z_l,z_r]$ of width $w_I$ and midpoint $m_I$, the mean-value theorem gives
\begin{equation}
\sup_{z\in I}|\varepsilon(z)|\le |\varepsilon(m_I)|+\frac{w_I}{2}\sup_{z\in I}|\varepsilon'(z)|.
\label{eq:cert_bound}
\end{equation}
Both terms on the right of Equation~\eqref{eq:cert_bound} are enclosed rigorously by ball arithmetic. Intervals that do not yet satisfy the target are bisected adaptively. The logarithmic endpoint $z=0$ is treated separately by an analytic mean-value bound, avoiding direct evaluation of $\ln z$ on an interval containing zero.

The exact derivative needed by the certifier follows from Equation~\eqref{eq:Gexact_hypergeom}:
\begin{equation}
\frac{dG_{\mathrm{exact}}}{dz}=-\frac{\pi}{4}\,{}_2F_1\!\left(\frac12,\frac12;2;1-z\right).
\label{eq:Gexact_derivative}
\end{equation}
This makes the verification independent of the exploratory minimax routine.

\subsection{Certified \texorpdfstring{$\ELIP\;R_{5,2}+3\log[T_4/F_3]$}{ELIP R52+3log}}
For the fixed coefficients in Equation~\eqref{eq:R52_cert_params}, the interval computation proves
\begin{equation}
1.6225893294455\times10^{-9}<\max|\varepsilon|<1.6229999624073\times10^{-9}.
\label{eq:R52_cert_rel}
\end{equation}
Equivalently,
\begin{equation}
0.0016225893294455<E_{\max}<0.0016229999624073\ppm,
\label{eq:R52_cert_ppm}
\end{equation}
so a concise certified statement is
\begin{equation}
\boxed{E_{\max}<0.0016230\ppm.}
\label{eq:R52_cert_box}
\end{equation}
The lower bound is supplied by a rigorous point evaluation near a numerically observed extremum; the upper bound covers the complete interval.

\subsection{Certified \texorpdfstring{$\ELIP\;R_{5,3}+3\log[T_8/F_0;\beta_1]$}{ELIP R53+3log beta1}}
For Equation~\eqref{eq:R53_beta1}, the exact quantity being certified is
\begin{equation}
G_{\beta}(h)=G_{53,0}(h)\left(1-Ah^{12}z^{7/4}\right),
\label{eq:R53_cert_function}
\end{equation}
with the decimal coefficients interpreted as fixed exact decimal rationals in the executable specification. A 192-bit Arb run with a separately treated interval near $z=0$ completed the adaptive cover and gave
\begin{equation}
4.51187500469109\times10^{-10}<\max|\varepsilon|<4.51999984378346\times10^{-10}.
\label{eq:R53_cert_rel}
\end{equation}
In ppm,
\begin{equation}
0.000451187500469<E_{\max}<0.000451999984379\ppm,
\label{eq:R53_cert_ppm}
\end{equation}
and therefore
\begin{equation}
\boxed{E_{\max}<0.000452\ppm.}
\label{eq:R53_cert_box}
\end{equation}
The successful run used 81,762 accepted subintervals after processing 163,012 nodes. A separate endpoint estimate on $0\le z\le10^{-14}$ was several orders of magnitude below the target. The certificate establishes the error bound of the fixed formula; it does not claim global minimax optimality of the fitted amplitude $A$.

\subsection{Relative standing of the certified high-accuracy formula}
Using the certified bound in Equation~\eqref{eq:R53_cert_box}, the maximum error is at least about $35$ times smaller than Popenkov's certified $0.015852\ppm$ bound, about $1.27\times10^3$ times smaller than the published R2/2exp numerical maximum, and about $3.20\times10^4$ times smaller than Cantrell's maximum. Relative to Ramanujan II the reduction is roughly $8.9\times10^5$. Against the particular ten-constant UNDERWOOD coefficient set evaluated in Table~\ref{tab:highcomparison}, the ELIP bound is about $1.44$ times lower numerically. These comparisons should not be interpreted as a universal ranking over all UNDERWOOD parameter budgets: the cited UNDERWOOD paper contains larger bases with substantially smaller reported errors.

Near the circle, the structural difference is stronger than the global numbers suggest. Because the $R_{5,3}$ core matches through $h^8$, the first nonzero base error is $O(h^9)$, while the external factor in Equation~\eqref{eq:R53_beta1} starts at $O(h^{12})$. Thus the correction does not consume the high-order circular information that was obtained analytically.

\section{Conclusions}\label{sec:conclusions}

This paper has developed the ELIP family from a deliberately minimal geometric starting point: the straight line joining the exact normalized circle value $\pi$ and the exact fully elongated value $4$. The central idea is to bend this line with correction structures whose ingredients are chosen to match known mathematics of the exact perimeter rather than to imitate the complete function blindly.

Two branches were explored. The logarithmic--power--exponential branch uses the correct qualitative branch-point behavior and then adds low-cost deformation or exponential modulators. The rational--logarithmic branch separates the correction into a rational smooth residual and exact logarithmic terms from the elongated-limit asymptotics. The notation $\ELIP\;R_{n,m}+q\log[T_r/F_f;\beta_t]$ makes explicit how much of a formula is analytically constrained and how much is genuinely fitted.

The low-complexity results show that useful accuracy does not require a large formula. The rationalized $R_{2,2}$ formula reaches about $3.57\ppm$ with arithmetic alone, the rationalized $R_{2,2}+1\log$ version reaches about $0.904\ppm$, and the purely rational $R_{3,3}$ construction with $b_3=29$ reaches about $1.596\ppm$ while retaining the same leading near-circle error as Cantrell. These formulas occupy a favorable simplicity--accuracy region and are independent of Ramanujan I and II.

The high-accuracy results illustrate a second principle: exact local information and global minimax compensation need not compete if they are assigned to different parts of the architecture. The certified $R_{5,2}+3\log[T_4/F_3]$ formula has $E_{\max}<0.001623\ppm$. When all coefficients of the rational core are instead fixed by exact Taylor matching, the error becomes smooth enough to be corrected efficiently by a beta-like factor. The $R_{5,3}+3\log[T_8/F_0;\beta_1]$ member retains a fully analytical rational core, uses a single fitted amplitude, and satisfies the rigorous global bound $E_{\max}<0.000452\ppm$.

Within the comparisons made here, this last formula is about $35$ times more accurate globally than Popenkov, roughly $1.27\times10^3$ times more accurate than the published R2/2exp formula, about $3.2\times10^4$ times more accurate than Cantrell, and nearly $8.9\times10^5$ times more accurate than Ramanujan II when maximum relative error is used as the metric. On the quoted ten-constant UNDERWOOD set it also has a smaller numerical global maximum, while the difference close to the circle is much larger because ELIP preserves exact Maclaurin information through $h^8$. The result should be interpreted carefully: larger UNDERWOOD bases reported in the literature reach still smaller errors, so the claim is not that ELIP dominates every closed-form construction at every parameter budget. It is also worth emphasizing an important structural difference between the two approaches. For instance, UNDERWOOD U10 employs ten coefficients obtained through numerical fitting, whereas $\mathrm{ELIP}\;R_{5,3}+3\log[T_8/F_0;\beta_1]$ contains a rational core whose eight coefficients are determined analytically from exact Taylor-series matching, with only a single parameter introduced through numerical minimax fitting. Moreover, as a consequence of its construction from the ELIP baseline, the latter approximation is exact at both limiting values of the shape parameter, $h=0$ and $h=1$. Thus, although the comparison of maximum errors is important, the two formulas also differ substantially in the amount of fitted information required and in the exact analytical constraints incorporated into their construction.

The main outcome is instead methodological. A highly accurate ellipse-perimeter approximation can be built without taking either Ramanujan formula as its base, by combining endpoint geometry, exact series coefficients, explicit logarithmic branch terms, and a very small number of minimax degrees of freedom. The certified $0.000452\ppm$ formula is therefore among the most accurate compact closed-form ellipse-perimeter approximations considered in the present comparison, while retaining an unusually transparent division between exact analytical information and numerical fitting.

Future work includes interval certification of additional ELIP members, symbolic simplification of the exact $\mathbb{Q}(\pi)$ coefficient sets, systematic degree-growth studies for fully constrained $R_{n,m}$ cores, and a theoretical analysis of the sign and unimodality properties of their beta-like residuals.


\section*{Funding}
This research received no external funding.

\section*{Data and Code Availability}
The formulas are fully specified in the manuscript. Numerical optimization, error evaluation,
figure generation, and rigorous interval-arithmetic certification scripts are available at \url{https://doi.org/10.5281/zenodo.23309365} so that all headline error bounds can be independently reproduced.

\section*{Declaration on the Use of Generative AI}

Generative AI tools were used during the preparation of this manuscript
to assist with LaTex formatting, exploratory numerical
coding, and the execution of some computational experiments. All
mathematical derivations, numerical results, references, figures, and
certified error bounds reported in the manuscript were reviewed and
verified by the authors, who take full responsibility for the content
of the work.

\section*{Conflicts of Interest}
The authors declare no conflicts of interest.

\bibliographystyle{unsrt}

\begin{thebibliography}{999}

\bibitem{AntennaElliptic}
Sivagnanam, S.; Kailasapathi, P.; Ramesh Kumar, S. A multi-band crescent shaped microstrip patch antenna with spiral shape ground slot for UWB applications. \emph{International Journal for Research in Engineering Application and Management} \textbf{2018}, \emph{4}, 489--493.

\bibitem{BrehlDowVAM}
Brehl, D. E.; Dow, T. A. Review of vibration-assisted machining. \emph{Precis. Eng.} \textbf{2008}, \emph{32}, 153--172.

\bibitem{CompositeEllipse}
Nakka, R.; Kumar, A. P.; Harursampath, D.; Ponnusami, S. A. Influence of fibre cross-section profile on the multi-physical properties of uni-directional composites. \emph{Composite Structures} \textbf{2023}, \emph{321}, 117321.

\bibitem{HumanThermalEllipse}
Sun, X. \emph{Development of an Improved Thermal Model of the Human Body and an Experimental Investigation of Heat Transfer from a Moving Cylinder}. Ph.D. Thesis, Kansas State University, Manhattan, KS, USA, 2012.

\bibitem{Kepler1609}
Kepler, J. \emph{Astronomia Nova}; Prague, 1609. Available online: \url{https://archive.org/details/Astronomianovaa00Kepl}.

\bibitem{Almkvist1988}
Almkvist, G.; Berndt, B. Gauss, Landen, Ramanujan, the arithmetic--geometric mean, ellipses, $\pi$, and the Ladies Diary. \emph{Am. Math. Mon.} \textbf{1988}, \emph{95}, 585--608.

\bibitem{Sykora2005}
Sýkora, J. Approximations of ellipse perimeters and of the complete elliptic integral $E(x)$: review of known formulae. \emph{Ebyte.it Stan's Library III}, 2005. Available online: \url{http://www.ebyte.it/library/docs/math05a/EllipsePerimeterApprox05.html}.

\bibitem{Kennedy1973}
Kennedy, H. C. \emph{Selected Works of Giuseppe Peano}; University of Toronto Press: Toronto, 1973.

\bibitem{Euler1774}
Euler, L. \emph{Nova Series Infinita Maxime Convergens Perimetrum Ellipsis Exprimens}; 1774. Available online: \url{https://scholarlycommons.pacific.edu/euler-works/448}.

\bibitem{Michon2007}
Michon, G. P. Final Answers: Perimeter of an Ellipse, 2007. Available online: \url{https://www.numericana.com/answer/ellipse.htm}.

\bibitem{Ramanujan1914a}
Ramanujan, S. Modular equations and approximations to $\pi$. \emph{Q. J. Math.} \textbf{1914}, \emph{45}, 350--372.

\bibitem{Koshy2024}
Koshy, K. I. Ellipse perimeter approximation: Two high accuracy formulae derived from a perimeter property. \emph{Int. J. Sci. Res. Math. Stat. Sci.} \textbf{2024}, \emph{11}, 1--5.

\bibitem{Moscato2024}
Moscato, P.; Ciezak, A. A new approximation for the perimeter of an ellipse. \emph{Algorithms} \textbf{2024}, \emph{17}, 464. doi:10.3390/a17100464.

\bibitem{Sykora2007Advances}
Sýkora, S. Advances in approximations of ellipse perimeters and of the complete elliptic integral. \emph{Stan's Library}, 2007.

\bibitem{Powell1981}
Powell, M. J. D. \emph{Approximation Theory and Methods}; Cambridge University Press: Cambridge, 1981.

\bibitem{Ayala2026}
Ayala-Raggi, S. E.; Rendón-Marín, M. An exponential correction to Ramanujan's second formula for ellipse perimeter computation. \emph{AppliedMath} \textbf{2026}, \emph{6}, 56. doi:10.3390/appliedmath6040056.

\bibitem{Popenkov2026}
Popenkov, O. A compact rational--logarithmic approximation to the perimeter of an ellipse. \emph{Zenodo}, 2026. doi:10.5281/zenodo.21798982.

\bibitem{Mamoun2_2026}
Mamoun, M. Two Foci, One Cup Runneth Over Again: The Ellipse Perimeter, Cracked Harder Than Ever---The Definitive Second Breakthrough in Closed-Form Approximator Models. \emph{Zenodo}, 2026. doi:10.5281/zenodo.19830852.

\end{thebibliography}

\end{document}